\documentclass[journal,twoside,print]{ieeecolor}
\usepackage{generic}
\usepackage{cite}
\usepackage{amsmath,amssymb,amsfonts}
\usepackage{graphicx}
\usepackage{algorithm,algorithmic}
\usepackage{hyperref}
\hypersetup{hidelinks=true}
\usepackage{textcomp}
\def\BibTeX{{\rm B\kern-.05em{\sc i\kern-.025em b}\kern-.08em
    T\kern-.1667em\lower.7ex\hbox{E}\kern-.125emX}}
\usepackage{subfig}
\usepackage{epstopdf}
\usepackage{multirow} 
\usepackage{threeparttable}
\usepackage{bbding}

\definecolor{gray}{rgb}{0.95, 0.95, 0.95}

\newtheorem{Thm}{Theorem}
\newtheorem{Def}{Definition}
\newtheorem{Lem}{Lemma}
\newtheorem{Ass}{Assumption}
\newtheorem{Rem}{Remark}
\newtheorem{Cor}{Corollary}

\begin{document}
\title{Pareto-optimal Trade-offs Between Communication and Computation
\\
with Flexible Gradient Tracking}

\author{Yan Huang$^\dagger$$^\ddagger$, Jinming Xu$^\dagger$, Li Chai$^\dagger$, \IEEEmembership{Senior Member, IEEE}, Jiming Chen$^\dagger$, \IEEEmembership{Fellow, IEEE}, 
\\
and Karl Henrik Johansson$^\ddagger$, \IEEEmembership{Fellow, IEEE}
\thanks{$^\dagger$College of Control Science
and Engineering, Zhejiang University, 310027 Hangzhou, China (e-mail: {jimmyxu@zju.edu.cn, chaili@zju.edu.cn, cjm@zju.edu.cn}).}
\thanks{$^\ddagger$Division of Decision
and Control Systems, School of EECS, KTH Royal Institute of Technology, SE-100 44 Stockholm, Sweden (e-mail: {yahuang@kth.se, kallej@kth.se}). }
}

\maketitle

\thispagestyle{headings}

\begin{abstract}
This paper addresses distributed stochastic optimization problems under non-i.i.d. data, focusing on the inherent trade-offs between communication and computational efficiency. 
To this end, we propose FlexGT, a flexible snapshot gradient tracking method that enables tunable numbers of local updates and neighbor communications per round, thereby adapting efficiently to diverse system resource conditions.
Leveraging a unified convergence analysis framework, we derive tight communication and computational complexity for FlexGT with explicit dependence on objective properties and certain tunable parameters.
Moreover, we introduce an accelerated variant, termed Acc-FlexGT, and prove that, with prior knowledge of the graph, it achieves Pareto-optimal trade-offs between communication and computation.
Particularly, in the nonconvex case, Acc-FlexGT achieves the optimal iteration complexity of $\tilde{\mathcal{O}}\left( \left( L\sigma ^2 \right) /\left( n\epsilon ^2 \right) +L/\left( \epsilon \sqrt{1-\sqrt{\rho _W}} \right) \right) $ and optimal communication complexity of $\tilde{\mathcal{O}}\left( L/\left( \epsilon \sqrt{1-\sqrt{\rho _W}} \right) \right)$ for appropriately chosen numbers of local updates, matching existing lower bounds up to logarithmic factors. And, it improves the existing results for the strongly convex case by a factor of $\tilde{\mathcal{O}} \left( 1/\sqrt{\epsilon} \right)$, where $\epsilon$ is the targeted accuracy, $n$ the number of nodes, $L$ the Lipschitz constant, $\rho_W$ the connectivity of the graph, and $\sigma$ the stochastic gradient variance. Numerical experiments corroborate the theoretical results and demonstrate the effectiveness of the proposed methods.
\end{abstract}


\section{Introduction}

With rapid growth in both data volume and device scale \cite{konevcny2016federated}, distributed optimization methods have gained increasing attention due to their wide range of applications
such as in cooperative control \cite{nedic2018distributed}, distributed sensing \cite{akyildiz2011cooperative}, large-scale machine learning \cite{tsianos2012consensus}, distributed training of large language models (LLMs) \cite{jia2024sdp4bit}, to name a few. In this paper, we consider a distributed stochastic optimization problem jointly solved by $n$ nodes over a network:
\begin{equation}\label{Prob} 
\underset{x\in \mathbb{R} ^p}{\min}f\left( x \right) =\frac{1}{n}\sum_{i=1}^n{\underset{:=f_i\left( x \right)}{\underbrace{\mathbb{E} _{\xi _i\sim \mathcal{D} _i}\left[ f_i\left( x;\xi _i \right) \right] }}},
\end{equation}
where $x\in \mathbb{R} ^p$ is the global decision variable. The objective function $f: \mathbb{R} ^p\rightarrow \mathbb{R}$ is the sum of $n$ local objective functions $f_i$, $i=1,\dots ,n$, conditioned on the local data sample $\xi _i$ with distribution $\mathcal{D} _i$.

\subsection{Literature Review}
Decentralized stochastic gradient descent (DSGD) \cite{ram2009asynchronous} is a simple yet efficient algorithm for solving Problem~\eqref{Prob}. It performs stochastic gradient descent locally while communicating with neighbors to achieve consensus, which enables data parallelism and avoids communication bottlenecks in centralized architectures \cite{mcmahan2017communication}. However, the native DSGD methods often struggle in scenarios with non-independent and identically distributed (non-i.i.d.) local datasets \cite{Li2020On, gorbunov2021local}. For Problem~\eqref{Prob}, such a scenario corresponds to the case where the local data distributions $\mathcal{D}_i$ differ across agents. 
This data heterogeneity induces gradient drift in the sense that local gradients computed at each node deviate from the global gradient direction, necessitating decaying stepsizes to ensure exact convergence \cite{nedic2009distributed, Li2020On}. This requirement limits the applicability of DSGD in many real-world tasks, as non-i.i.d. data distributions are common.

Many efforts have been made to address the issue of data heterogeneity in distributed optimization.
Gradient tracking (GT) methods have gained popularity for mitigating the effect \cite{xu2015augmented, di2016next, shi2015extra, li2020revisiting, sun2022distributed, tian2022acceleration}. Notably, Xu \textit{et al.} \cite{xu2015augmented} proposed AugDGM, which introduces an auxiliary variable at each node to track the average of the local gradients, enabling exact convergence with a constant stepsize. This method has been extended to broader settings, including time-varying graphs \cite{nedic2017achieving} and directed graphs \cite{pu2020push, saadatniaki2020decentralized}, with convergence guarantees. Pu \textit{et al.} \cite{pu2021distributed} extended GT to the stochastic setting, introducing distributed stochastic gradient tracking (DSGT). They proved that DSGT achieves linear convergence to a neighborhood of the optimal solution, determined by the stochastic gradient noise. This extension enhances the applicability of GT, for instance, to many machine learning tasks. Although GT methods are robust to data heterogeneity, they exhibit a stronger dependence on the graph compared to DSGD \cite{huang2022tackling, koloskova2021improved}, leading to increased communication overhead.
The consensus alternating direction method of multipliers (C-ADMM) \cite{ling2015dlm} and its inexact variant \cite{chang2014multi} also exhibit robustness to data heterogeneity. However, they require either solving a subproblem exactly at each iteration, which incurs high computation overhead, or assuming a well-structured objective function to enable linear approximations, which limits their applicability.
These observations highlight the need for more efficient algorithms, especially in resource-constrained scenarios that involve both communication and computation.

To mitigate the high communication burden, many improved GT-based methods have been developed \cite{nguyen2023performance, song2022compressed, song2024optimal, didouble, wu2025effectiveness, liu2025decentralized}.
A popular strategy is to skip some communication steps and perform local updates only. For instance, FedAvg \cite{mcmahan2017communication} (also known as Local SGD \cite{stich2018local}) is widely used in federated learning \cite{dean2012large, karimireddy2020scaffold}, where periodic communication occurs between server and worker nodes, with partial node participation to reduce communication overhead.
In decentralized settings, Nguyen \textit{et al.} \cite{nguyen2023performance} introduced LU-GT, a variant of GT that uses multiple local updates, and provided communication complexity results matching those of Local SGD for nonconvex objectives. 
The effectiveness of multiple local updates under mild data heterogeneity is further demonstrated in \cite{wu2025effectiveness}.
Subsequently, Liu \textit{et al.} \cite{liu2025decentralized} proposed the gradient-sum-tracking method K-GT, which reduces communication complexity by achieving lower stochastic gradient variance. However, both \cite{nguyen2023performance} and \cite{liu2025decentralized} mainly focus on communication efficiency and do not explicitly consider how increasing the number of local updates may raise computational complexity due to weakened inter-node consensus.
Another approach is to randomly or periodically communicate over a subset of the graph’s edges to mitigate the impact of single-point communication bottleneck and stragglers \cite{ram2009asynchronous, lian2017can, lian2018asynchronous, ying2021exponential}.
These methods can be seen as approximate asynchronous update schemes, which often require specific connectivity assumptions \cite{ram2009asynchronous} or particular topologies, such as exponential graphs \cite{ying2021exponential}.

In contrast to directly skipping certain communications, as in Local SGD \cite{stich2018local}, another strategy involves using accelerated consensus protocols. For instance, Chebyshev acceleration \cite{song2024optimal, didouble} and accelerated gossip \cite{liu2011accelerated, lu2021optimal, yuan2022revisiting} methods enable the distributed algorithm to behave more like a centralized algorithm. In particular, Song \textit{et al.} \cite{song2024optimal} introduced an optimal gradient tracking (OGT) method with Chebyshev acceleration for deterministic optimization problems with smooth and strongly convex objective functions. Leveraging prior knowledge of the network topology, they demonstrated that a limited number of communication steps per round could enable the distributed algorithm to replicate centralized updates. Building on this, the stochastic variant SS-DSGT \cite{didouble} was proposed, achieving lower communication complexity compared to DSGT. For nonconvex cases, Lu \textit{et al.} \cite{lu2021optimal} developed an accelerated DSGT method called DeTAG, achieving optimal \textit{iteration} complexity with large batches of stochastic gradients. 
Although communication quantization and compression can be used to further reduce the amount of data transmitted in each round \cite{koloskova2019decentralized, song2022compressed}, their impacts on algorithmic accuracy must be carefully considered. For a comprehensive overview of other efforts to improve communication efficiency, readers are referred to the recent survey papers \cite{tang2020communication, cao2023communication}.

Some recent efforts have focused on both the communication and computational efficiency of distributed optimization algorithms.
Berahas {\textit{et al.}} \cite{berahas2018balancing} proposed a distributed gradient descent method with multiple communication steps in each round, named NEAR-DGD, and evaluated the performance of the algorithm using a weighted combination of the number of communication and computation steps.
In \cite{berahas2024balancing}, the authors explored the use of multiple local updates and inter-node communication in each round with a deterministic GT method, showing that increasing the number of communication steps improves the convergence rate for the strongly convex objective functions.
However, in distributed stochastic optimization, balancing communication and computation is more challenging, as stochastic gradient noise critically affects the convergence rate and renders results from deterministic methods inapplicable \cite{berahas2018balancing, berahas2024balancing}. Although the results in \cite{iakovidou2022s, liu2022decentralized} highlight the impact of communication and computation protocols on \textit{iteration} complexity, the trade-off between these two factors is not explicitly characterized, particularly in the presence of data heterogeneity. Additionally, most of the studies focus on specific types of objective functions, such as strongly convex or nonconvex objectives, or rely on strict assumptions, such as bounded gradients \cite{liu2022decentralized}. We summarize and compare the most relevant algorithms in Table~\ref{Tab_comparison_of_related_work}.

Part of the results in this work were previously presented at a conference \cite{huang2023computation}. The current work: i) extends the algorithm in \cite{huang2023computation} by incorporating snapshot GT and accelerated gossip communication, resulting in a more flexible framework with improved convergence guarantees; ii) develops a unified theoretical framework for strongly convex, convex, and nonconvex objective functions with comprehensive complexity results; and iii) includes more extensive numerical experiments on both synthetic and real-world datasets.

\begin{table*}
    \begin{center}
        \caption{Relevant algorithms for solving Problem~\eqref{Prob}}
        \label{Tab_comparison_of_related_work}
        {
        \begin{threeparttable}
            \begin{tabular}{c|c|c|c|c|c|c}
            \hline
                \rule{0pt}{10pt}
                {\textbf{Reference}}
                &{\textbf{Objectives}}
                & {\textbf{Communication}}
                & {\textbf{Computation}}
                &{\textbf{Stochastic}}
                & {\textbf{Acceleration}}
                & {\textbf{Non-i.i.d.}}
                \\
                \hline
                DSGT\cite{pu2021distributed, koloskova2021improved} & SC, Cvx, NC  & Single & Single & {\checkmark} & {\tiny \XSolidBrush} & {\checkmark}
                \rule{0pt}{10pt}
                \\
                \hline
                \rule{0pt}{10pt}
                 SS-GT \cite{song2024optimal}&  SC  & Single & Single &  {\tiny \XSolidBrush} & {\checkmark} & {\checkmark}
                 \\
                \hline
                \rule{0pt}{10pt}
                 SS-DSGT \cite{didouble}&  SC & Single & Single & {\checkmark} & {\checkmark} & {\checkmark}
                \\
                \hline
                \rule{0pt}{10pt}
                LU-GT\cite{nguyen2023performance}&  NC  & Single & Multiple & {\tiny \XSolidBrush} & {\tiny \XSolidBrush} & {\checkmark}
                \\
                \hline
                \rule{0pt}{10pt}
                K-GT \cite{liu2025decentralized}&  NC  & Single & Multiple & {\checkmark} & {\tiny \XSolidBrush} & {\checkmark}
                \\
                \hline
                \rule{0pt}{10pt}
                DeTAG \cite{lu2021optimal}&  NC  & Multiple & Single  & {\checkmark} & {\checkmark} & {\checkmark}
                \\
                \hline
                \rule{0pt}{10pt}
                GTA \cite{berahas2024balancing} & SC & Multiple& Multiple 
                & {\tiny \XSolidBrush} &  {\tiny \XSolidBrush} & {\checkmark}
                 \\
                \hline
                \rule{0pt}{10pt}
                 S-NEAR-DGD\cite{iakovidou2022s} & SC& Multiple & Multiple & {\checkmark} & {\tiny \XSolidBrush} & {\tiny \XSolidBrush} 
                 \\
                \hline
                \rule{0pt}{10pt}
                 DFL \cite{liu2022decentralized} & NC& Multiple & Multiple & {\checkmark} & {\tiny \XSolidBrush} & {\tiny \XSolidBrush} 
                 \\
                \hline
                \rule{0pt}{10pt}
                (Acc-)FlexGT (\textbf{ours})& SC, Cvx, NC& Multiple & Multiple & {\checkmark} & {\checkmark} & {\checkmark}
                \\
                \hline
            \end{tabular}
            
            \begin{tablenotes}
               \item Note: SC, Cvx, and NC denote strongly convex, convex, and nonconvex objective functions, respectively. “Single” and “Multiple” indicate that the algorithm uses single or multiple local communication or computation steps per round, respectively. “Stochastic” refers to the use of stochastic gradients. “Acceleration” refers to communication acceleration.
            \end{tablenotes}
        \end{threeparttable}
        }
    \end{center}
\end{table*}

\subsection{Contributions}

This paper solves the distributed stochastic optimization problem \eqref{Prob} in scenarios with non-i.i.d local datasets. We introduce the FlexGT algorithm and illustrate the inherent trade-offs between communication and computation. Pareto-optimal trade-offs (cf., Definition \ref{Def_Pareto_optimum}) are further achieved through an accelerated variant, termed Acc-FlexGT. The contributions are summarized as follows:
\begin{itemize}
    \item We propose FlexGT, a flexible gradient tracking method that allows adjustable numbers of local updates and communications in each round, enabling better adaptation to heterogeneous communication and computation resources in distributed systems. By introducing a snapshot point within the GT step, FlexGT mitigates the inconsistency between GT variables during local updates, leading to improved convergence. Moreover, by integrating an accelerated gossip communication protocol, we develop Acc-FlexGT, which significantly enhances both communication and computational efficiency when the graph connectivity parameter is known a priori. To the best of our knowledge, Acc-FlexGT is the first algorithm to achieve Pareto-optimal trade-offs between communication and computation from a multi-objective optimization perspective, jointly minimizing these two costs under given resource budgets.
    \item We develop a unified convergence analysis framework for FlexGT and Acc-FlexGT, covering objective functions that are strongly convex, convex, or nonconvex. 
    By employing refined error decomposition and decoupling techniques, we establish linear or sublinear convergence rates with the best-known communication and computational complexity in most settings, explicitly characterizing the dependence on problem-specific parameters, compared to existing GT algorithms (Theorem~\ref{Thm_FlexGT} and Corollary~\ref{Col_complexity_FlexGT}). In particular, we show that Acc-FlexGT achieves the optimal iteration complexity and optimal communication complexity in the nonconvex case with appropriately chosen numbers of local updates, matching existing lower bounds up to logarithmic factors (Corollary~\ref{Cor_Acc-FlexGT}). Moreover, Acc-FlexGT improves the best-known results in the strongly convex case by a factor of $\tilde{\mathcal{O}}(1/\sqrt{\epsilon})$. This unified analysis framework can be extended to a broader class of algorithms, which is of independent interest.
\end{itemize}

\subsection{Paper Organization}
The remainder of this paper is structured as follows. Section~\ref{Sec_prob} formulates the distributed stochastic optimization problem over networks under standard assumptions and introduces the FlexGT and Acc-FlexGT algorithms. Section~\ref{Sec_main_results} presents the main convergence result of the algorithms for strongly convex, convex, and nonconvex objective functions, showing the dependency on both the amount of communication and computation. Section~\ref{Sec_conv_analysis} provides the proof of the result. Numerical experiments on synthetic and real-world datasets are given in Section~\ref{Sec_experiments} to validate the effectiveness of FlexGT and Acc-FlexGT. Finally, Section~\ref{Sec_conclusion} concludes the paper.

\subsection{Notations} Throughout this paper, we adopt the following notations: $\left\| \cdot \right\|$ represents the Frobenius norm, $\left< \cdot ,\cdot \right> $  the inner product of vectors, $\mathbb{E}\left[ \cdot \right]$  the expectation of a random variable, $\mathbf{1}$ the all-ones vector, $\mathbf{I}$ the identity matrix, and $\mathbf{J}=\mathbf{1}\mathbf{1}^{\top}/n$ the averaging matrix. 
Moreover, we use the asymptotic notation $\mathcal{O}(\cdot)$ to hide constant factors, $\tilde{\mathcal{O}}(\cdot)$ to additionally suppress logarithmic terms, and $\lceil \cdot \rceil$ to represent the ceiling function.

\section{Problem Formulation and Algorithm Design}\label{Sec_prob}
To solve the distributed stochastic optimization problem~\eqref{Prob} in a decentralized manner, we consider the following equivalent problem:
\begin{equation}\label{Prob_cons}
\begin{aligned}
\underset{\mathbf{x}\in \mathbb{R} ^{n\times p}}{\min}&F\left( \mathbf{x} \right) =\frac{1}{n}\sum_{i=1}^n{\underset{:=f_i\left( x_i \right)}{\underbrace{\mathbb{E} _{\xi _i\sim \mathcal{D} _i}\left[ f_i\left( x_i;\xi _i \right) \right] }},}
\\
&\mathrm{s.t.}\,\,x_i=x_j,\quad \forall i,j\in \left[ n \right] 
\end{aligned}
\end{equation}
where $\mathbf{x}:=\left[ x_1,x_2,\dots, x_n \right]^\top \in \mathbb{R} ^{n\times p}$ denotes the collection of the local copies of the decision variables from all nodes. The constraints guarantee that problems \eqref{Prob} and \eqref{Prob_cons} are equivalent.
Moreover, the agents are connected over a graph $\mathcal{G}=(\mathcal{V},\mathcal{E})$. Here, $\mathcal{V}=\{1,2,...,n\}$ represents the set of nodes, and $\mathcal{E}\subseteq\mathcal{V}\times \mathcal{V}$ denotes the set of edges consisting of ordered pairs $(i,j)$ representing the communication link from node $j$ to node $i$. For node $i\in[n]$, we define $\mathcal{N}_i= \{ j| \left( i, j \right) \in \mathcal{E} \}$ as the set of its neighboring nodes (including $i$ itself).

\begin{figure}[t]
    \centering
    \subfloat 
    {
        \begin{minipage}[t]{0.45\textwidth}
            \centering
            \includegraphics[width=\textwidth]{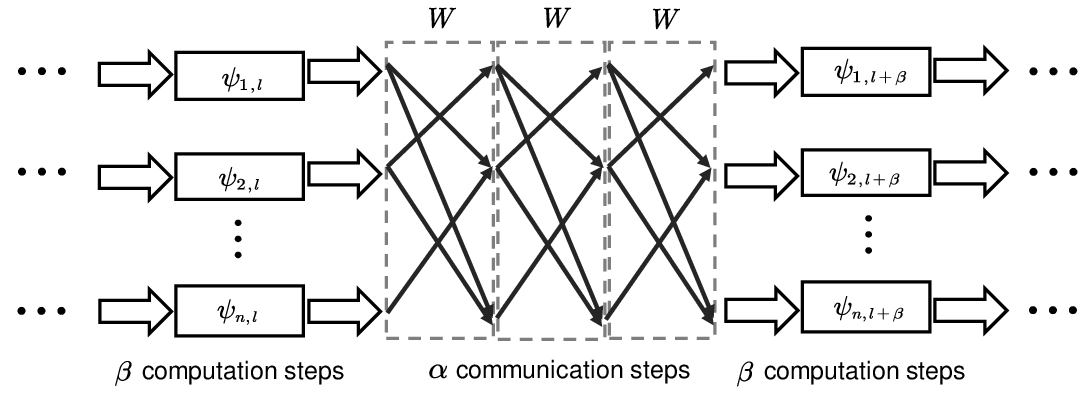}
        \end{minipage}%
    }
    \caption{Flow diagram of the computation and communication of FlexGT. The variable $\psi_{i,l}$, $i=1,\dots,n$, represents the collection of $x_{i,l}$, $y_{i,l}$, and $z_{i,l}$. Only $x_{i,l}$ and $y_{i,l}$ are communicated. The communication is weighted by the matrix $W$.}
    \label{Fig_flowchart}
\end{figure}

\subsection{Assumptions and Definitions}
We make the following assumptions regarding the objective function $f_i$, its gradient $\nabla f_i$, and the graph $\mathcal{G}$.

\begin{Ass}[\textbf{Convexity}]\label{Ass_convexity}
For each $f_i$, there is a constant $\mu \geqslant 0$ such that for all  $x, x^{\prime}\in \mathbb{R} ^p$,
\begin{equation*}
f_i\left( x^{\prime} \right) -f_i\left( x \right) \geqslant \left< \nabla f\left( x \right) ,x^{\prime}-x \right> +\frac{\mu}{2}\left\| x^{\prime}-x \right\| ^2.
\end{equation*}
\end{Ass}

Particularly, we say that $f_i$ is $\mu$-strongly convex if $\mu > 0$; otherwise, we refer to it as being convex.
\begin{Ass}[\textbf{Smoothness}]\label{Ass_smooth}
Each $f_i $ is $L$-smooth, i.e., there is a Lipschitz constant $L>0$ such that for all  $x, x^{\prime} \in \mathbb{R} ^p$,
    \begin{equation*}
        f_i\left( x^{\prime} \right) -f_i\left( x \right) \leqslant \left< \nabla f\left( x \right) ,x^{\prime}-x \right> +\frac{L}{2}\left\| x-x^{\prime} \right\| ^2.
    \end{equation*}
Moreover, the overall objective function is bounded from below, i.e., $f\left( x \right) \geqslant f^*:=\mathrm{inf}_{x\in \mathbb{R} ^p}\,\,f\left( x \right) >-\infty$.
\end{Ass}

In our algorithms, we assume that each agent $i$ obtains an unbiased noisy gradient of the form $\nabla f_i\left( x; \xi _i \right) $ by querying a stochastic oracle ($\mathcal{S}\mathcal{O}$) with $\xi _i\sim \mathcal{D} _i$ at each iteration, where $\mathcal{D} _i$ is the local data distribution. This is typically achieved by uniformly random sampling on each node.

\begin{Ass}[\textbf{Bounded variance}]\label{Ass_bounded_var}
For the stochastic gradient $\nabla f_i\left( x; \xi _i \right) $ with $\xi _i\sim \mathcal{D} _i$, there exist a constant $\sigma  \geqslant 0$ such that for all $x\in \mathbb{R}^p$,
\begin{equation*}
\begin{aligned}
\mathbb{E} \left[ \left\| \nabla f_i\left( x;\xi _i \right) -\nabla f_i\left( x \right) \right\| ^2 \right] \leqslant \sigma ^2.
\end{aligned}
\end{equation*}
\end{Ass}

\begin{Ass}[\textbf{Graph connectivity}] \label{Ass_graph}
The weight matrix $W$ induced by the graph $\mathcal{G}$ is doubly stochastic, i.e., $W\mathbf{1}=\mathbf{1},\mathbf{1}^{\top}W=\mathbf{1}^{\top}$, and $\rho _W:= \left\| W-\mathbf{J} \right\|_2 ^2 <1.$
\end{Ass}

Note that a suitable weight matrix  $W$  satisfying Assumption~\ref{Ass_graph} can be constructed for any undirected and connected graph. For example, a weight matrix can be readily determined using the Metropolis-Hastings protocol \cite{xiao2006distributed}.

In this work, we aim to jointly minimize communication and computational complexities, measured respectively by the total number of inter-node communications and local stochastic gradient evaluations needed to achieve a prescribed accuracy. This naturally leads to a multi-objective optimization problem that admits a set of Pareto-optimal solutions, as defined below.

\begin{Def}[\textbf{Pareto-optimality}~\cite{boyd2004convex, chen2013distributed}]\label{Def_Pareto_optimum}
For a multi-objective optimization problem that minimizes $\left\{ g_j\left( \omega \right) \right\} _{j\in \left[ m \right]}$, a point $\omega$ is Pareto-optimal if there does not exist another point $\omega ^{\prime}$ such that: $\forall j \in [m]$, $g_j\left( \omega ^{\prime} \right) \leqslant g_j\left( \omega \right)$; and $\exists j\in \left[ m \right] $, $g_j\left( \omega ^{\prime} \right) < g_j\left( \omega \right)$. The set of Pareto-optimal points is called the Pareto frontier.
\end{Def}

In what follows, we refer to the Pareto-optimal points that jointly optimize communication and computational complexity in terms of tunable parameters as Pareto-optimal trade-offs. 

It should be noted that the above multi-objective optimization perspective is fundamentally different from the existing \emph{iteration-complexity–based} analyses \cite{lu2021optimal, yuan2022revisiting}, which characterizes the minimum number of algorithmic iterations required to reach a target accuracy, under a given function class, algorithm class, and stochastic oracle model.
Such settings do not explicitly account for communication and computational costs per iteration, both of which are, indeed, influenced by the use of multiple communication rounds, multiple local updates, or large-batch stochastic gradients.

\subsection{FlexGT and Acc-FlexGT Algorithms}
Now, we proceed to present the FlexGT algorithm and its accelerated variant, Acc-FlexGT, for solving problem \eqref{Prob_cons}. The pseudo-code for the algorithms is provided in Algorithm~\ref{Alg_Flexible_DSGT}. Note that the integers $\alpha, \beta \geqslant 1$ denote the number of communication and computation steps performed in each round, respectively. A flow diagram illustrating the communication and computation of FlexGT is shown in Fig.~\ref{Fig_flowchart}.

\begin{algorithm}[t]
\caption{\textbf{FlexGT and Acc-FlexGT}}
\begin{algorithmic}[1]
\label{Alg_Flexible_DSGT}
\renewcommand{\algorithmicrequire}{\textbf{Initialization: }}
\REQUIRE {Initial points $z_{i,0}=x_{i,0}\in \mathbb{R}^p$ and $y_{i,0}=\nabla f_i\left( x_{i,0};\xi _{i,0} \right)$, number of communication and computation steps $\alpha , \beta \geqslant 1$, and stepsize $\gamma >0$.
}
\
\FOR{round $k = 0,1,\dots$, each node $i\in[n]$,}

\STATE Update snapshot point: $z _{i,\beta k}=x_{i,\beta k}$.

\FOR{$l=\beta k, \beta k+1,\dots , \beta \left( k+1 \right) -1$}
\STATE Update decision variable:
\begin{equation*}
    x_{i,l+1}=x_{i,l}-\gamma y_{i,l}.
\vspace{-0.5cm}
\end{equation*}
\STATE Let $z_{i,l+1}=z_{i,l}$.
\STATE Compute the stochastic gradient\\ $g_{i,l+1}=\nabla f_i\left( z_{i,l+1};\xi _{i,l+1} \right)$ by querying $\mathcal{S}\mathcal{O}$.
\STATE Update gradient tracking variable: 
\[y_{i,l+1}=y_{i,l} + g_{i,l+1} - g_{i,l}.\]
\ENDFOR

\STATE Let $\theta_{i,0}=\theta_{i,-1}=\left\{ x_{i,\beta \left( k+1 \right)}, y_{i,\beta \left( k+1 \right)} \right\} $.

\FOR{$s = 0,1,\dots,\alpha -1$}
\STATE \textbf{FlexGT}: Direct inter-node communication:

\begin{equation*}
\theta_{i,s+1}= \sum_{j\in \mathcal{N} _i}{w_{ij}\theta_{j,s}}.
\end{equation*}
\STATE \textbf{Acc-FlexGT}: Accelerated gossip communication:
\begin{equation*}
\theta_{i,s+1}=\left( 1+\eta \right) \sum_{j\in \mathcal{N} _i}{w_{ij}\theta_{j,s}-\eta}\theta_{i,s-1},
\end{equation*}
where $\eta =\frac{1-\sqrt{1-\rho _W}}{1+\sqrt{1-\rho _W}}$.
\ENDFOR
\ENDFOR
\end{algorithmic}
\end{algorithm}

We first introduce the notations for the collections of decision variables, GT variables, and snapshot points, respectively:
\begin{equation*}
\begin{aligned}
&\mathbf{x}_t:=\left[ x_{1,t},x_{2,t},\dots ,x_{n,t} \right] ^{\top}\in \mathbb{R} ^{n\times p},
\\
&\mathbf{y}_t:=\left[ y_{1,t},y_{2,t},\dots ,y_{n,t} \right] ^{\top}\in \mathbb{R} ^{n\times p},
\\
&\mathbf{z}_t:=\left[ z_{1,t},z_{2,t},\dots ,z_{n,t} \right] ^{\top}\in \mathbb{R} ^{n\times p},
\end{aligned}
\end{equation*}
Then, Algorithm~\ref{Alg_Flexible_DSGT} can be rewritten in a compact form:
\begin{subequations}\label{Compact_PerDSGT}
\begin{align}
\mathbf{x}_{\beta \left( k+1 \right)}&=\bar{W}\left( \mathbf{x}_{\beta k}-\gamma \sum_{j=0}^{\beta -1}{\mathbf{y}_{\beta k+j}} \right) ,\label{Eq_x_update}
\\
\mathbf{y}_{\beta \left( k+1 \right)}&=\bar{W}\left( \mathbf{y}_{\beta k}+G_{\beta \left( k+1 \right)}-G_{\beta k} \right),\label{Eq_y_update}
\end{align}
\end{subequations}
where 
\[
G_t:=\left[ \dots ,g_{i,t},\dots \right] ^{\top}\in \mathbb{R} ^{n\times p}, \quad g_{i,t}=\nabla f_i\left( z_{i,t};\xi _{i,t} \right)
\]
denotes the collection of the stochastic gradient at the snapshot point $\mathbf{z}_t$ at iteration $t$, and its expectation is denoted by
\[
\nabla F_t:=\mathbb{E} \left[ G_t \right] =\left[ \dots ,\nabla f_i\left( z_{i,t} \right) ,\dots \right] ^{\top}\in \mathbb{R} ^{n\times p}.
\]
We note that the form of $\bar{W}$ varies with the communication protocol used, as specified in lines 11 and 12 of Algorithm~\ref{Alg_Flexible_DSGT}. In the case of FlexGT (cf., line 11), we have $\bar{W}=W^{\alpha }$. Otherwise, when the accelerated gossip communication protocol is used as in Acc-FlexGT (cf., line 12), we have $\bar{W}=M_{\alpha }$ with $M_{\alpha }$ obtained through $\alpha$ communication steps:
\begin{equation*}
\begin{aligned}
M_{s+1}=\left( 1+\eta \right) WM_s-\eta M_{s-1},\quad s \in \{0,1,\dots, \alpha -1  \}
\end{aligned}
\end{equation*}
with $M_{-1}=M_0=I$. Note that $\bar{W}$ is still a doubly-stochastic matrix. Although both cases involve $\alpha$ communications per round, they differ in the available prior knowledge, i.e., whether the spectral gap of the graph is known a priori for designing the stepsize $\eta$, which results in distinct communication properties, as stated in the following lemma.

\begin{Lem}\label{Lem_communication_protocol}
Suppose Assumption~\ref{Ass_graph} holds and denote $\bar{\rho}_W:=\left\| \bar{W}-\mathbf{J} \right\| ^2$. Then, for FlexGT, we have 
\begin{equation}\label{Eq_multi_comm}
\bar{\rho}_W=\left\| W^{\alpha }-\mathbf{J} \right\| ^2\leqslant \rho _{W}^{\alpha };
\end{equation}
and for Acc-FlexGT,
\begin{equation}\label{Eq_fastmix}
\bar{\rho}_W=\left\| M_{\alpha }-\mathbf{J} \right\| ^2\leqslant 2\left( 1-\sqrt{1-\sqrt{\rho _W}} \right) ^{2\alpha }.
\end{equation}
\end{Lem}

Equation \eqref{Eq_multi_comm} is straightforward to derive based on Assumption~\ref{Ass_graph}, and \eqref{Eq_fastmix} was stated at Proposition 3 in \cite{liu2011accelerated}. The difference in $\bar{\rho}_W$ reflects the difference in the consensus speed achieved by different communication strategies under the same number of communication rounds.

Regarding the snapshot point $\mathbf{z}$, we have
\begin{equation}\label{Eq_z_update}
{\mathbf{z}_{\beta k+t}}=\,\,\mathbf{x}_{\beta k}, \quad t=0,1,\dots ,\beta -1,
\end{equation}
which implies that $\mathbf{y}$ tracks the stochastic gradient evaluated at the snapshot point $\mathbf{z}$, rather than at the current decision variable $\mathbf{x}$ as in conventional GT methods \cite{pu2021distributed}, i.e., 
\begin{equation}
\mathbf{1}^{\top}{\mathbf{y}_{\beta k+t}}=\mathbf{1}^{\top}{G_{\beta k+t}}=\sum_{i=1}^n{\nabla f_i\left( x_{i,\beta k};\xi _{i,\beta k+t} \right)}.
\end{equation}
This prevents the tracking variables from diverging excessively due to multiple local updates, i.e., the gradient difference $G_{t+1}-G_t$ remains controlled, which is crucial for stability in the presence of heterogeneous data. 
This method differs from the existing snapshot GT approach \cite{song2024optimal, didouble}, which randomly skips updating the tracking variables while performing their consensus at each iteration. In contrast, FlexGT utilizes the stochastic gradient at each iteration to slightly update the tracking variables and perform communication periodically to avoid extra costs from redundant initialization.

The FlexGT algorithm generalizes many existing methods by incorporating an adjustable communication and computation protocol parameterized by $\alpha$ and $\beta$. Specifically, it recovers DSGT \cite{pu2021distributed} ($\alpha = \beta = 1$), LU-GT \cite{nguyen2023performance} ($\alpha = 1, \beta \geqslant 1$), and DeTAG ($\alpha \geqslant 1, \beta = 1$) as special cases. It can also recover DFL \cite{liu2022decentralized} when an independent weight matrix is used for updating $\mathbf{y}$. This highlights the flexibility of FlexGT and its applicability to a broader range of problems.

\section{Main Results}\label{Sec_main_results}
In this section, we present the main convergence results of the FlexGT and Acc-FlexGT algorithms for strongly convex, convex, and nonconvex objective functions under a unified analysis framework. 
To this end, we first define the following Lyapunov function for (strongly) convex cases, consisting of the optimality gap, consensus error, and GT error:
\begin{equation}\label{Def_Lyapunov_func}
\begin{aligned}
V_t:=\left\| \bar{x}_t-x^* \right\| ^2+c_x\left\| \tilde{\mathbf{x}}_t \right\| ^2+c_y\left\| \tilde{\mathbf{y}}_t \right\| ^2,
\end{aligned}
\end{equation}
where $x^{*}$ indicates the global optimal solution of Problem~\eqref{Prob}, $t$ denotes the iteration, $c_x$ and $c_y$ are coefficients to be properly designed later, and the consensus error and GT error are defined as follows:
\begin{equation*}
\begin{aligned}
\tilde{\mathbf{x}}_{t}&:=\mathbf{x}_t-\mathbf{1}\bar{x}_{t}, \quad \bar{x}_t:=\mathbf{1}^{\top}\mathbf{x}_t/n,
\\
\tilde{\mathbf{y}}_t&:=\mathbf{y}_t-\mathbf{1}\bar{y}_t, \quad \bar{y}_t:=\mathbf{1}^{\top}\mathbf{y}_t/n.
\end{aligned}
\end{equation*}

\begin{table*}
    \footnotesize 
    \begin{center}
        \caption{Communication and computational complexity of FlexGT required to reach an accuracy level of $\epsilon > 0$ for strongly convex, convex, and nonconvex objective functions.}
        \label{Tab_complexity_FlexGT}
        \resizebox{1\textwidth}{!}{
        \begin{threeparttable}
            \begin{tabular}{c|c|c|c}
                \hline
                \rule{0pt}{10pt}
                {\textbf{Algorithm}}&
                {\textbf{Objective}} & {\textbf{Computation} ($\tilde{\mathcal{O}} \left( \cdot \right) $)} & {\textbf{Communication} ($\tilde{\mathcal{O}} \left( \cdot \right) $) }
                \\
                \hline
                \rule{0pt}{15pt}
                &
                \multirow{3}{*}{SC} & \multirow{3}{*}{$
                   \frac{\beta L}{\left( 1-\rho _{W}^{\alpha } \right) ^2\mu}+\frac{\sigma ^2}{\mu ^2n\epsilon}+\frac{\beta \sqrt{L\rho _{W}^{\alpha }\sigma ^2}}{\mu ^{3/2}\left( 1-\rho _{W}^{\alpha } \right) ^{3/2}\epsilon ^{1/2}}$} 
                & 
                \multirow{3}{*}{$\frac{\alpha L}{\left( 1-\rho _{W}^{\alpha } \right) ^2\mu}+\frac{\alpha \sigma ^2}{\beta \mu ^2n\epsilon}+\frac{\alpha \sqrt{L\rho _{W}^{\alpha }\sigma ^2}}{\mu ^{3/2}\left( 1-\rho _{W}^{\alpha } \right) ^{3/2}\epsilon ^{1/2}}$}           
                \\  & & &
                \\  & & &
                \rule{0pt}{10pt}
                \\
                \rule{0pt}{10pt}
                \multirow{3}{*}{FlexGT}
                &
                \multirow{3}{*}{Cvx} & \multirow{3}{*}{$\frac{\beta L}{\left( 1-\rho _{W}^{\alpha } \right) ^2\epsilon}+\frac{\sigma ^2}{n\epsilon ^2}+\frac{\beta \sqrt{L\rho _{W}^{\alpha }\sigma ^2}}{\left( 1-\rho _{W}^{\alpha } \right) ^{3/2}\epsilon ^{3/2}}$} 
                & 
                \multirow{3}{*}{$\frac{\alpha L}{\left( 1-\rho _{W}^{\alpha } \right) ^2\epsilon}+\frac{\alpha \sigma ^2}{n\beta \epsilon ^2}+\frac{\alpha \sqrt{L\rho _{W}^{\alpha }\sigma ^2}}{\left( 1-\rho _{W}^{\alpha } \right) ^{3/2}\epsilon ^{3/2}}$}
                \\  & & &
                \\ & & &
                \rule{0pt}{10pt}
                \\
                \rule{0pt}{10pt}
                &
                \multirow{3}{*}{NC}
                &
                \multirow{3}{*}{$\frac{L\sigma ^2}{n\epsilon ^2}+\frac{\beta L}{\left( 1-\rho _{W}^{\alpha } \right) ^2\epsilon}+\frac{\beta R_0}{\epsilon}+\frac{\beta L\sqrt{\rho _{W}^{\alpha }\sigma ^2}}{\left( 1-\rho _{W}^{\alpha } \right) ^{3/2}\epsilon ^{3/2}}$}
                & 
                \multirow{3}{*}{$\frac{\alpha L\sigma ^2}{n\beta \epsilon ^2}+\frac{\alpha L}{\left( 1-\rho _{W}^{\alpha } \right) ^2\epsilon}+\frac{\alpha R_0}{\epsilon}+\frac{\alpha L\sqrt{\rho _{W}^{\alpha }\sigma ^2}}{\left( 1-\rho _{W}^{\alpha } \right) ^{3/2}\epsilon ^{3/2}}$}                 
                \\  & & &
                \\  & & &
                \rule{0pt}{15pt}
                \\
                \hline
               
            \end{tabular}
        \end{threeparttable}
        }
    \end{center}
\end{table*}

\begin{table*}
    \footnotesize 
    \begin{center}
        \caption{Communication and computational complexity of Acc-FlexGT required to reach an accuracy level of $\epsilon > 0$ for strongly convex, convex, and nonconvex objective functions. The last column ``Iteration$^*$" presents the corresponding best-known or optimal \textit{iteration} complexity.}
        \label{Tab_complexity_Acc_FlexGT}
        \resizebox{1\textwidth}{!}{
        \begin{threeparttable}
            \begin{tabular}{c|c|c|c|c}
            \hline
                \rule{0pt}{10pt}
                {\textbf{Algorithm}}&
                {\textbf{Objective}} & {\textbf{Computation} ($\tilde{\mathcal{O}} \left( \cdot \right) $)} & {\textbf{Communication} ($\tilde{\mathcal{O}} \left( \cdot \right) $) } & {\textbf{Iteration$^*$} ($\tilde{\mathcal{O}}\left( \cdot \right) $ )}
                \\
                \hline
                \rule{0pt}{15pt}
                &
                \multirow{2}{*}{SC} & \multirow{2}{*}{$
                    \frac{\sigma ^2}{n\mu ^2\epsilon}+ \frac{\beta L}{\mu}$} 
                & 
                \multirow{2}{*}{$\left( \frac{\sigma ^2}{n\beta \mu ^2\epsilon}+\frac{L}{\mu} \right) \frac{1}{\sqrt{1-\sqrt{\rho _W}}}$}  &    \multirow{2}{*}{$\frac{\sigma ^2}{n\mu ^2\epsilon}+\frac{\sqrt{L\sigma ^2}}{\mu \sqrt{1-\sqrt{\rho _W}}\epsilon ^{1/2}}$ \cite{didouble} \tnote{a}}   
                \\  & & & &
                \rule{0pt}{10pt}
                \\
                \rule{0pt}{10pt}
                \multirow{2}{*}{Acc-}
                &
                \multirow{2}{*}{Cvx} & \multirow{2}{*}{$\frac{L\sigma ^2}{n\epsilon ^2}+\frac{\beta \sqrt{L}}{\epsilon}$} 
                &
                \multirow{2}{*}{$\left( \frac{L\sigma ^2}{n\beta \epsilon ^2}+\frac{\sqrt{L}}{\epsilon} \right) \frac{1}{\sqrt{1-\sqrt{\rho _W}}}$} & \multirow{2}{*}{/ \tnote{b}}
                
                \\ \text{FlexGT}
                & & & & 
                \rule{0pt}{10pt}
                \\
                \rule{0pt}{10pt}
                &
                \multirow{2}{*}{NC}
                &
                \multirow{2}{*}{$\frac{L\sigma ^2}{n\epsilon ^2}+\frac{\beta L}{\epsilon}$}
                &
                \multirow{2}{*}{$\left( \frac{L\sigma ^2}{n\beta \epsilon ^2}+\frac{L}{\epsilon} \right) \frac{1}{\sqrt{1-\sqrt{\rho _W}}}$}  & \multirow{2}{*}{$\frac{L\sigma ^2}{n\epsilon ^2}+\frac{L}{\sqrt{1-\sqrt{\rho _W}}\epsilon}$ \cite{lu2021optimal} \tnote{c} }              
                \\  & & & & 
                \rule{0pt}{15pt}
                \\
                \hline
            \end{tabular}
            \begin{tablenotes}
                \item[a] To the best of our knowledge, this is the best-known result in the existing literature in this setting.
                \item[b] The mark "/" indicates no relevant result in the existing literature.
                \item[c] This is the optimal iteration complexity that matches the existing lower bound up to a logarithmic factor.
            \end{tablenotes}
        \end{threeparttable}
        }
    \end{center}
\end{table*}

\subsection{Convergence}
The following theorem presents the convergence results of FlexGT and Acc-FlexGT for strongly convex, convex, and nonconvex objective functions under specific assumptions.

\begin{Thm}\label{Thm_FlexGT}
Suppose Assumptions~\ref{Ass_smooth}--\ref{Ass_graph} hold. There exists a constant stepsize 
\begin{equation}\label{Eq_stepsize_global}
\gamma =\mathcal{O} \left( \min \left\{ 1, \frac{\left( 1-\bar{\rho}_W \right) ^2}{\bar{\rho}_W} \right\} \right) \frac{1}{\beta L}
\end{equation}
such that the optimization error is bounded as follows:
\\
\textbf{i) Strongly convex:} Under the additional Assumption~\ref{Ass_convexity} with $\mu > 0$,  we have for all $k\geqslant0$
\begin{equation}\label{Eq_convergence_FlexGT_sc}
\begin{aligned}
\mathbb{E} \left[ V_{\beta \left( k+1 \right)} \right] &\leqslant \underset{\mathrm{linear \,\,rate}}{\underbrace{\left( 1-\min \left\{ \frac{\mu \beta \gamma}{2},\frac{1-\bar{\rho}_W}{8} \right\} \right) }}\mathbb{E} \left[ V_{\beta k} \right] 
\\
&\quad +\underset{\mathrm{centralized \,\, error}}{\underbrace{\frac{\gamma ^2\beta \sigma ^2}{n}}}+\underset{\mathrm{topology \,\,effect}}{\underbrace{\frac{1664\gamma ^3\beta ^{3}L\bar{\rho}_W}{\left( 1-\bar{\rho}_W \right) ^3}\sigma ^2}}.
\end{aligned}
\end{equation}
\textbf{ii) Convex:} Under the additional Assumption~\ref{Ass_convexity} with $\mu \geqslant 0$, we have for all $K\geqslant1$,
\begin{equation}\label{Eq_convergence_FlexGT_gc}
\begin{aligned}
&\frac{1}{K}\sum_{k=0}^{K-1}{\mathbb{E} \left[ f\left( \bar{x}_{\beta k} \right) -f\left( x^* \right) \right]}
\\
&\leqslant \underset{\mathrm{sublinear\,\, rate}}{\underbrace{\frac{2\mathbb{E} \left[ V_0 \right]}{\gamma \beta K}}}+\underset{\mathrm{centralized \,\, error} \mathrm{}}{\underbrace{\frac{2\gamma \sigma ^2}{n}}}
+\underset{\mathrm{topology \,\, effect}}{\underbrace{\frac{3328\gamma ^2\beta ^{2}L\bar{\rho}_W}{\left( 1-\bar{\rho}_W \right) ^3}\sigma ^2}}.
\end{aligned}
\end{equation}
\textbf{iii) Nonconvex:} It holds that for all $K\geqslant1$,
\begin{equation}\label{Eq_convergence_FlexGT_nc}
\begin{aligned}
&\frac{1}{K}\sum_{k=0}^{K-1}{\mathbb{E} \left[ \left\| \nabla f\left( \bar{x}_{\beta k} \right) \right\| ^2 \right] }
\\
&\leqslant\underset{\mathrm{sublinear\,\,rate}}{\underbrace{\frac{8\mathbb{E} \left[ f\left( \bar{x}_0 \right) -f^* \right]}{\gamma \beta K}+\frac{8\gamma ^2\beta ^{2}L^2\bar{\rho}_W\mathbb{E} \left[ \left\| \tilde{\mathbf{y}}_0 \right\| ^2 \right]}{n\left( 1-\bar{\rho}_W \right) ^3K}}}
\\
&\quad+\underset{\mathrm{centralized \,\, error}}{\underbrace{\frac{4\gamma L \sigma ^2}{n}}}+\underset{\mathrm{topology \,\, effect}}{\underbrace{\frac{3328\gamma ^2\beta ^{2}L^2\bar{\rho}_W}{\left( 1-\bar{\rho}_W \right) ^3}\sigma ^2}}.
\end{aligned}
\end{equation}
\end{Thm}

\begin{proof}
    See Section ~\ref{Sec_conv_analysis}.
\end{proof}

Theorem~\ref{Thm_FlexGT} shows that for the strongly convex case, with constant stepsize $\gamma$, FlexGT converges linearly to a neighborhood of the optimal solution of Problem~\eqref{Prob}. The steady-state error consists of two parts: one is the term of stochastic gradient noise matching that of the centralized SGD algorithm, while the other part is due to the decentralized graph topology. 
Similarly, in the convex case with $\mu=0$, the rate is sublinear in the sense of the running sum of the objective function. Finally, in the nonconvex case, the running sum of the gradient norm of FlexGT converges to a neighborhood of a stationary point at a sublinear rate. Note that these results highlight a clear dependence on the parameters $\mu$ and $L$ of the objective function, the computation steps $\beta$, and the spectrum gap $\bar{\rho}_W$ related to communication steps $\alpha$ and $\rho_W$ in Assumption~\ref{Ass_graph}.

\begin{Rem}\label{Rem_best_known_results}
Compared to the existing results in \cite{liu2022decentralized, iakovidou2022s}, FlexGT handles data heterogeneity \cite{lian2017can}, i.e. $\forall x$, 
\[
\zeta_f := \sup \left\{ \frac{1}{n} \sum_{i=1}^{n} \left\| \nabla f_i(x) - \nabla f(x) \right\|^2 \right\},
\]
and does not require the assumption of uniformly bounded stochastic gradients \cite{liu2022decentralized}, making it more robust in non-i.i.d. settings. Compared to the best-known results for DSGT ($\alpha = \beta = 1$) without acceleration in \cite{koloskova2021improved}, our findings recover their results in all three settings \footnote{Under Assumption~\ref{Ass_graph} on the graph and using the same convergence metric as in this work.} and extend to more general cases with $\alpha \geqslant 1$ and $\beta \geqslant 1$.
Compared K-GT \cite{liu2025decentralized} in the nonconvex setting, the communication complexity obtained in this work improves upon that of K-GT by a factor of $1/\sqrt{1-\rho _W} $ in terms of the topology effect, whereas K-GT benefits from a lower level of stochastic gradient noise.
\end{Rem}

Based on the convergence results in Theorem~\ref{Thm_FlexGT}, we can further derive the iteration complexity of FlexGT for different objective functions in the following corollary.

\begin{Cor}\label{Col_complexity_FlexGT}
Suppose Assumptions~\ref{Ass_smooth}--\ref{Ass_graph} hold. Let the stepsize satisfy the condition \eqref{Eq_stepsize_global} given in Theorem~\ref{Thm_FlexGT}, and denote by $\hat{K}$ the number of rounds required for FlexGT to reach an arbitrary accuracy $\epsilon > 0$. Then, 
\\
\textbf{i) Strongly convex:} Under the additional Assumption~\ref{Ass_convexity} with $\mu > 0$,  we have
\begin{equation}\label{Eq_complexity_FlexGT_sc}
\begin{aligned}
\hat{K}={\mathcal{O}}\Biggl( &\frac{L}{\left( 1-\bar{\rho}_W \right) ^2 \mu}\log \frac{\mathbb{E}\left[ V_0 \right]}{\epsilon}+\frac{\sigma ^2}{\beta \mu ^2n\epsilon} 
\\
&+\frac{\sqrt{L\bar{\rho}_W\sigma ^2}}{\mu ^{3/2} \left( 1-\bar{\rho}_W \right) ^{3/2}\epsilon ^{1/2}} \Biggl) .
\end{aligned}
\end{equation}
\textbf{ii) Convex:} Under the additional Assumption~\ref{Ass_convexity} with $\mu \geqslant 0$, we have
\begin{equation}\label{Eq_complexity_FlexGT_gc}
\begin{aligned}
\hat{K}=\mathcal{O} \left( \frac{L}{\left( 1-\bar{\rho}_W \right) ^2 \epsilon}+\frac{\sigma ^2}{\beta n\epsilon ^2}+\frac{\sqrt{L\bar{\rho}_W\sigma ^2}}{\left( 1-\bar{\rho}_W \right) ^{3/2}\epsilon ^{3/2}} \right) \mathbb{E} \left[ V_0 \right] .
\end{aligned}
\end{equation}
\textbf{iii) Nonconvex:} It holds
\begin{equation}\label{Eq_complexity_FlexGT_nc}
\begin{aligned}
\hat{K}=\mathcal{O}\Biggl(& \frac{L}{\left( 1-\bar{\rho}_W \right) ^2 \epsilon}+\frac{R_0}{\epsilon} +\frac{L\sigma ^2}{\beta n\epsilon ^2} 
\\
&+\frac{L\sqrt{\bar{\rho}_W\sigma ^2}}{\left( 1-\bar{\rho}_W \right) ^{3/2}\epsilon ^{3/2}} \Biggl) \left( \mathbb{E} \left[ f\left( \bar{x}_0 \right) -f^* \right] \right) ,
\end{aligned}
\end{equation}
where $R_0:=\frac{1}{n}\mathbb{E} \left[ \left\| \tilde{\mathbf{y}}_0 \right\| ^2 \right] \left( \mathbb{E} \left[ f\left( \bar{x}_0 \right) -f^* \right] \right) ^{-1}$.
\end{Cor}

\begin{proof}
See Section \ref{Proof_Cor_1}
\end{proof}

Based on the iteration complexity derived in Corollary~\ref{Col_complexity_FlexGT} and the number of local updates and communications per iteration, we can derive the communication and computational complexity of FlexGT for different objective functions, as summarized in Table~\ref{Tab_complexity_FlexGT}. To the best of our knowledge, FlexGT achieves the best-known complexity among GT-based methods for solving Problem~\eqref{Prob} in most settings (cf. Remark~\ref{Rem_best_known_results}).

\begin{Rem}\label{Rem_trade-off}
Table~\ref{Tab_complexity_FlexGT} highlights the impact of $\alpha$ and $\beta$ on the communication and computational complexity of FlexGT. Unlike previous works \cite{iakovidou2022s, liu2022decentralized}, our results reveal an inherent trade-off between communication and computation.
Particularly, as the number of computation steps $\beta$ increases, the communication complexity decreases while the computational complexity increases. As the number of communication steps $\alpha$ increases, the computational complexity decreases, while the communication complexity may either increase or decrease depending on $\rho_W$. This highlights the need for a careful design of $\alpha$ and $\beta$ to balance communication and computational costs according to the specific application requirements.
\end{Rem}

\subsection{Pareto-optimal Trade-offs}

\begin{figure}[!tb]
    \centering
    \subfloat[strongly convex] 
    {
        \begin{minipage}[t]{0.24\textwidth}
            \centering
            \includegraphics[width=\textwidth]{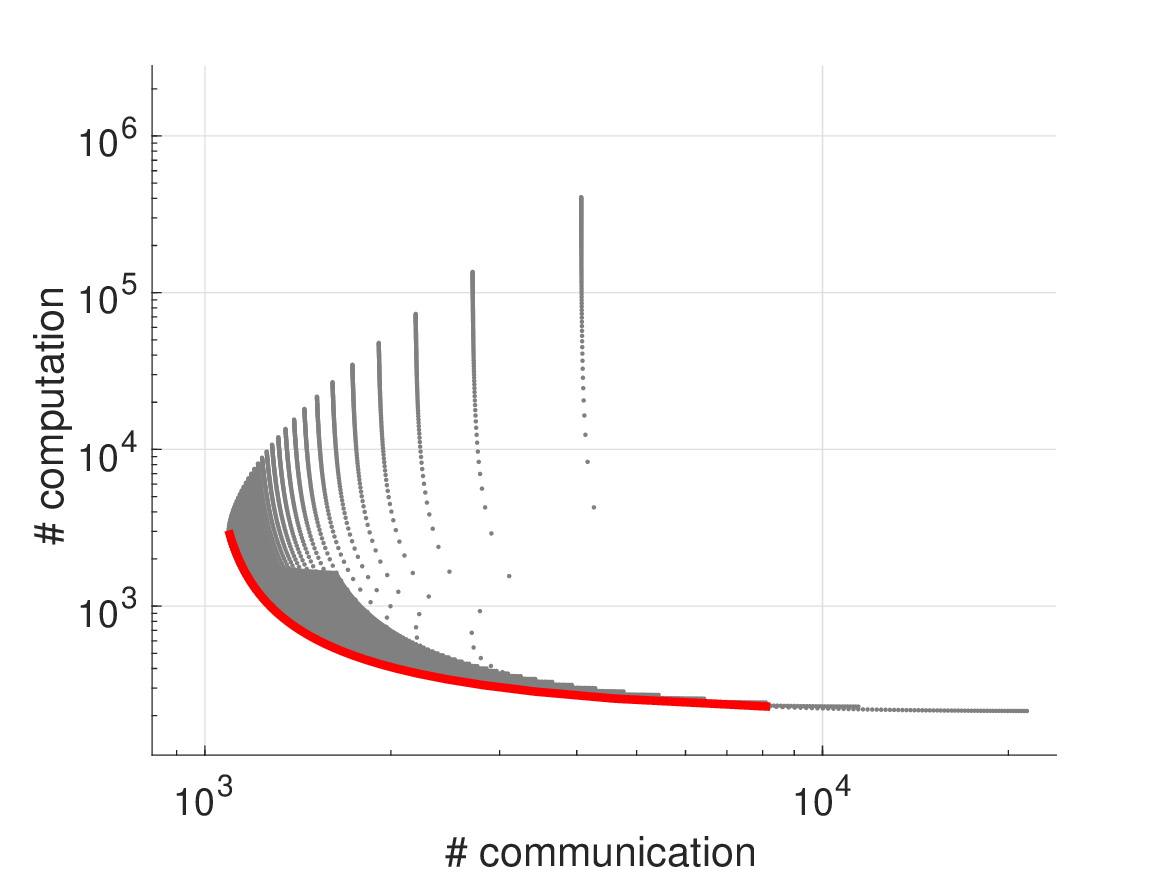} 
        \end{minipage}
    }
    \subfloat[nonconvex]
    {
        \begin{minipage}[t]{0.24\textwidth}
            \centering
            \includegraphics[width=\textwidth]{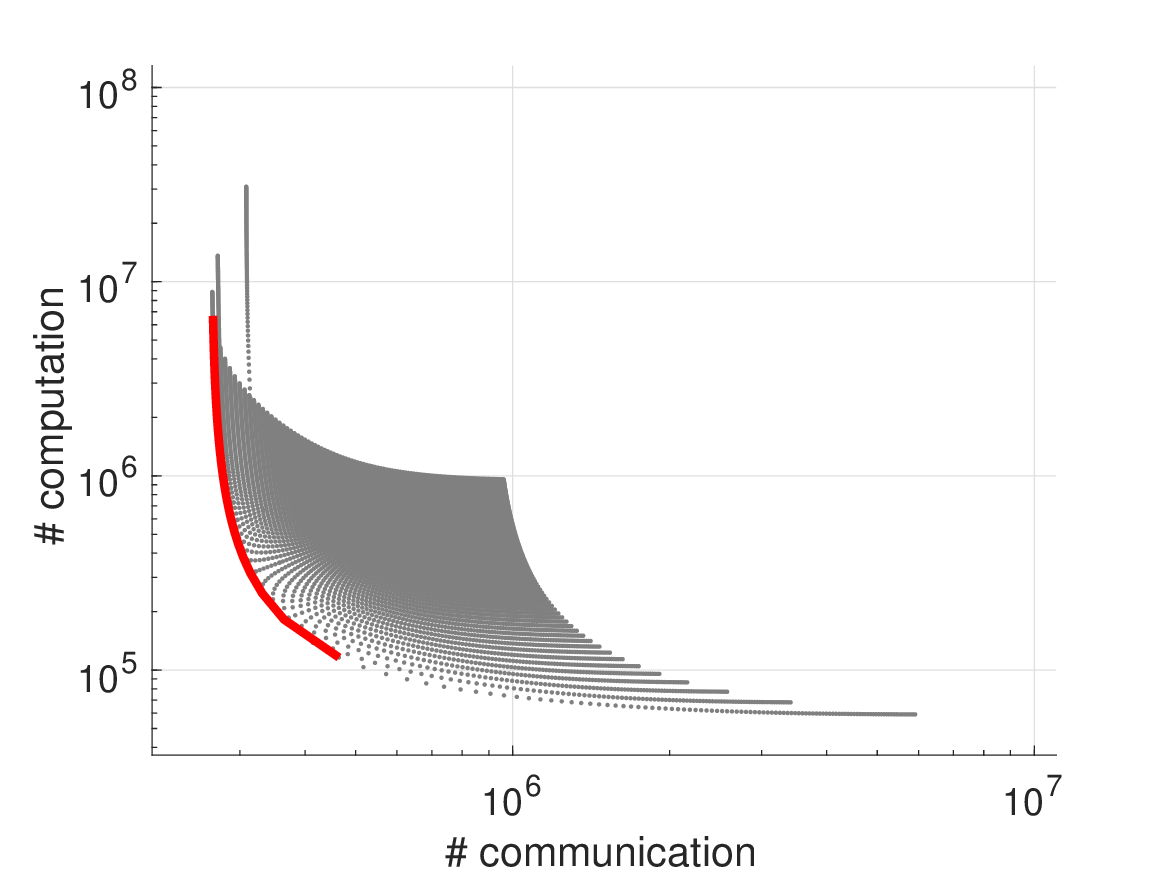}
        \end{minipage}
    }   
    \caption{The number of communication and computation steps needed for FlexGT to achieve an accuracy of $\epsilon = 10^{-4}$ with $\alpha, \beta =1,2,\dots,100$. Pareto-optimal solutions (red lines) are achieved with $\alpha = 4$ and $\alpha = 32$ for strongly convex (left) and nonconvex (right) cases, respectively.}
    \label{Fig_Pareto}
\end{figure}

The trade-off between communication and computation is illustrated in Fig.~\ref{Fig_Pareto}, where each gray point reflects a specific choice of $\alpha$ and $\beta$ in the range from 1 to 100. The coordinates of each point represent the estimated number of communication steps (horizontal axis) and number of computation steps (vertical axis) required to achieve an accuracy of $\epsilon$, obtained from Table~\ref{Tab_complexity_FlexGT}. It can be observed that, under Definition \ref{Def_Pareto_optimum} of Pareto-optimality in multi-objective optimization, a specific set of parameters lies on the Pareto frontier, as shown by the bold red lines. This set achieves the lowest overall complexity for reaching a given accuracy level, meaning that improvements in either communication or computational complexity are not possible without compromising the other.

The following corollary shows that by incorporating accelerated gossip communication, Acc-FlexGT substantially improves the iteration complexity and achieves Pareto-optimal trade-offs between communication and computation.

\begin{Cor}\label{Cor_Acc-FlexGT}
Suppose Assumptions~\ref{Ass_smooth}--\ref{Ass_graph} hold. Let the stepsize satisfy the condition \eqref{Eq_stepsize_global} given in Theorem~\ref{Thm_FlexGT}, and denote by $\hat{K}$ the number of rounds required for Acc-FlexGT to reach an arbitrary accuracy $\epsilon > 0$. Then,
\\
\textbf{i) Strongly convex:} Under the additional Assumption~\ref{Ass_convexity} with $\mu > 0$,  letting
\begin{equation}\label{Eq_d1_Acc_FlexGT_sc}
\alpha =\lceil \frac{\max \left\{ \ln 2,1/2\ln \left( n\beta L/\mu \right) \right\}}{\sqrt{1-\sqrt{\rho _W}}} \rceil,
\end{equation}
we have
\begin{equation}\label{Eq_complexity_Acc_FlexGT_sc}
\begin{aligned}
\hat{K}=\mathcal{O} \left( \frac{L}{\mu}\log \frac{\left[ V_0 \right]}{\epsilon}+\frac{\sigma ^2}{n\mu ^2\beta \epsilon} \right). 
\end{aligned}
\end{equation}
\textbf{ii) Convex:} Under the additional Assumption~\ref{Ass_convexity} with $\mu \geqslant 0$, letting
\begin{equation}\label{Eq_d1_Acc_FlexGT_gc}
\alpha =\lceil\frac{\max \left\{ \ln 2,1/2\ln \left( n\beta  \right) \right\}}{\sqrt{1-\sqrt{\rho _W}}}\rceil,
\end{equation}
we have
\begin{equation}\label{Eq_complexity_Acc_FlexGT_gc}
\begin{aligned}
\hat{K}=\mathcal{O} \left( \frac{\sqrt{L}}{\epsilon}+\frac{L\sigma ^2}{n\beta \epsilon ^2} \right) \mathbb{E} \left[ V_0 \right].
\end{aligned}
\end{equation}
\textbf{iii) Nonconvex:} Letting
\begin{equation}\label{Eq_d1_Acc_FlexGT_nc}
\alpha =\lceil{\frac{\max \left\{ \ln 2,1/2\ln \left( \beta \max \left\{ n,R_0 \right\} \right) \right\}}{\sqrt{1-\sqrt{\rho _W}}}}\rceil
\end{equation}
we have
\begin{equation}\label{Eq_complexity_Acc_FlexGT_nc}
\begin{aligned}
\hat{K}=\mathcal{O} \left( \frac{L}{\epsilon}+\frac{L\sigma ^2}{n\beta \epsilon ^2} \right) \mathbb{E} \left[ f\left( \bar{x}_0 \right) -f^* \right] .
\end{aligned}
\end{equation}
\end{Cor}

\begin{proof}
See Section \ref{Proof_cor_2}
\end{proof}

Based on Corollary~\ref{Cor_Acc-FlexGT}, the communication and computational complexities of Acc-FlexGT are obtained by multiplying $\hat{K}$ by $\alpha$ and $\beta$, respectively, as summarized in Table~\ref{Tab_complexity_Acc_FlexGT}. Compared to FlexGT, Acc-FlexGT achieves Pareto-optimal trade-offs with respect to $\beta$ by properly selecting $\alpha$, in the sense that it attains the Pareto frontier of the multi-objective optimization problem that jointly minimizes communication and computational costs.
In other words, any reduction in communication complexity entails an increase in computational complexity, and vice versa, as indicated by the red curves in Fig.~\ref{Fig_Pareto}.
This distinguishes our work from existing methods with accelerated communication \cite{lu2021optimal, yuan2022revisiting, didouble}, which achieve optimal \emph{iteration complexity} but do not explicitly account for the communication and computational costs incurred per iteration. Instead, we provide a more flexible framework to achieve Pareto-optimal trade-offs between these two costs via the tunable parameters $\alpha$ and $\beta$, and recover the results of the aforementioned methods as special cases, as detailed below.

\begin{Rem}\label{Rem_optimality}
Compared to existing results, for the nonconvex case, the obtained computational complexity of Acc-FlexGT
\begin{equation}
\tilde{\mathcal{O}}\left( \frac{L\sigma ^2}{n\epsilon ^2}+\frac{L}{\sqrt{1-\sqrt{\rho _W}}\epsilon} \right)
\end{equation}
is optimal in the sense that it matches the lower bound on the iteration complexity in \cite{lu2021optimal, yuan2022revisiting}, up to logarithmic factors arising from the choice of 
\begin{equation}
\beta =\alpha =\tilde{\mathcal{O}}\left( \frac{1}{\sqrt{1-\sqrt{\rho _W}}} \right).
\end{equation}
Moreover, Acc-FlexGT achieves the optimal communication complexity
\begin{equation}
\tilde{\mathcal{O}}\left( \frac{L}{\sqrt{1-\sqrt{\rho _W}}\epsilon } \right),
\end{equation}
with a sufficient number of local updates, i.e., 
\begin{equation}
\beta \geqslant \lceil \frac{\sigma ^2}{n\epsilon} \rceil,
\end{equation}
matching the lower bounds established for the deterministic case in \cite{sun2019distributed, yuan2022revisiting}.
In the strongly convex case, our result improves the first term in \eqref{Eq_complexity_Acc_FlexGT_sc} from $1/\sqrt{\epsilon}$ to $\log \left( 1/\epsilon \right) $ compared to \cite{didouble}, as we use a constant stepsize instead of a stepsize that decays with iterations. To the best of our knowledge, in the convex case, this work is the first to establish results for GT methods with local updates and Chebyshev acceleration.
\end{Rem}

\section{Convergence Analysis}\label{Sec_conv_analysis}
This section presents a unified convergence analysis framework for the theoretical results stated in Theorem~\ref{Thm_FlexGT}.

\begin{figure}[t]
    \centering
    \subfloat 
    {
        \begin{minipage}[t]{0.48\textwidth}
            \centering
            \includegraphics[width=\textwidth]{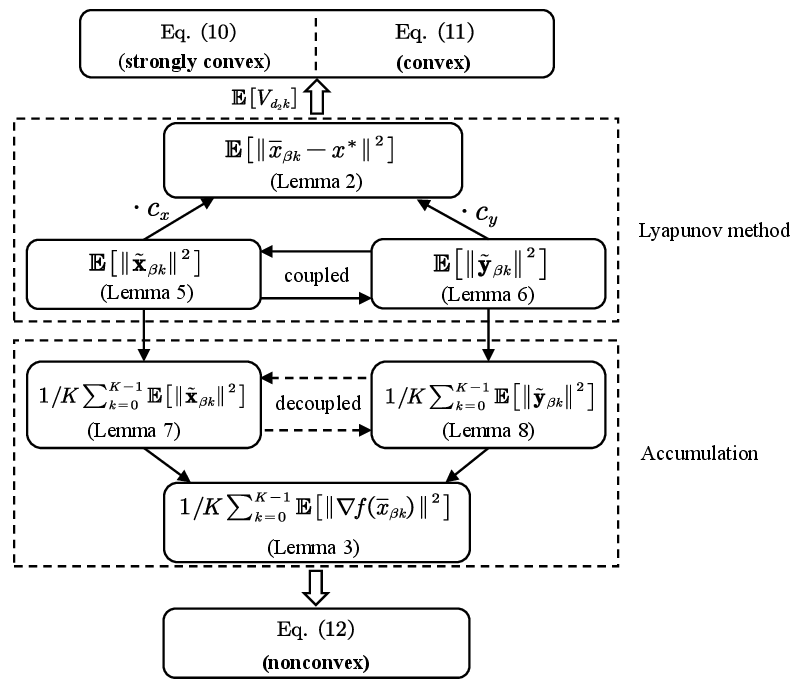}
        \end{minipage}%
    }
    \caption{Basic idea of the unified convergence analysis framework. The solid arrows show the dependency of the terms. The dashed arrows indicate decoupling between the two accumulated terms. The dashed box outlines the Lyapunov and the accumulation methods for (strongly) convex and nonconvex objective functions, respectively.}
    \label{Fig_proof}
\end{figure}

\textbf{Proof sketch.} The basic idea of the proof is illustrated in Fig.~\ref{Fig_proof}, which builds on the established contraction properties of the optimality gap for (strongly) convex objective functions (Lemma~\ref{Lem_opt_gap}) and descent property for the nonconvex case (Lemma~\ref{Lem_descent_FlexGT}), both of which involve dependencies on the consensus and GT errors. 
Rather than using a direct recursive analysis, the key to bounding these errors is to break the cyclic dependency between the consensus and GT errors (Lemmas~\ref{Lem_cons_err} and~\ref{Lem_tra_err}) as highlighted by the arrows in Fig.~\ref{Fig_proof}. For the (strongly) convex case, we achieve this by designing a suitable Lyapunov function in \eqref{Def_Lyapunov_func}. Conversely, for the nonconvex case, we accomplish this through a straightforward method of accumulating the errors (Lemmas~\ref{Lem_cons_err_sum_FlexGT} and~\ref{Lem_trac_err_sum_FlexGT}).
Notably, the proof relies on a double-loop analysis. We first derive the upper bounds of the client divergence over the interval between two communications (inner loop) (cf., Lemma~\ref{Lem_tech_client_diver} in the appendix). These bounds are subsequently utilized to establish the contraction properties of the consensus and GT errors during each communication round (outer loop).
Following this proof sketch, we begin by introducing some key lemmas that are crucial for the proof of Theorem~\ref{Thm_FlexGT}.

\subsection{Key Lemmas}

For the (strongly) convex case, i.e., when Assumption \ref{Ass_convexity} holds, we obtain the following lemma for the optimality gap.

\begin{Lem}\label{Lem_opt_gap}
Suppose Assumptions~\ref{Ass_convexity}--\ref{Ass_graph} hold. Let the stepsize satisfy $\gamma \leqslant \frac{1}{4\beta L}$. Then, we have for all $k \geqslant 0$,
\begin{equation}
\begin{aligned}
&\mathbb{E} \left[ \left\| \bar{x}_{\beta \left( k+1 \right)}-x^* \right\| ^2 \right] 
\\
&\leqslant \left( 1-\frac{\mu \beta \gamma}{2} \right) \mathbb{E} \left[ \left\| \bar{x}_{\beta k}-x^* \right\| ^2 \right] +\frac{\gamma ^2\beta \sigma ^2}{n}
\\
&\quad +\frac{3\gamma \beta L}{n}\mathbb{E} \left[ \left\| \tilde{\mathbf{x}}_{\beta k} \right\| ^2 \right] -\gamma \beta \mathbb{E} \left[ f\left( \bar{x}_{\beta k} \right) -f\left( x^* \right) \right] .
\end{aligned}
\end{equation}
\end{Lem}

\begin{proof}
According to the update rule of decision variable $\mathbf{x}$ \eqref{Eq_x_update}, we have
\begin{equation*}
\bar{x}_{\beta \left( k+1 \right)}=\bar{x}_{\beta k}-\gamma \sum_{j=0}^{\beta -1}{\frac{\mathbf{1}^{\top}}{n}\mathbf{y}_{\beta k+j}}.
\end{equation*}
Noticing that $\mathbf{1}^{\top}\bar{W}=\mathbf{1}^{\top}$ and the fact that $\mathbb{E} \left[ G_{\beta k} \right] =\nabla F_{\beta k}$, we get
\begin{equation}\label{Proof_lem1_1}
\begin{aligned}
&\mathbb{E} \left[ \left\| \bar{x}_{\beta \left( k+1 \right)}-x^* \right\| ^2 \right] 
\\
&=\mathbb{E} \left[ \left\| \bar{x}_{\beta k}-\gamma \sum_{j=0}^{\beta -1}{\frac{\mathbf{1}^{\top}}{n}\mathbf{y}_{\beta k+j}}-x^* \right\| ^2 \right] 
\\
&=\mathbb{E} \left[ \left\| \bar{x}_{\beta k}-x^* \right\| ^2 \right]
+\gamma ^2\underset{=:S_2}{\underbrace{\mathbb{E} \left[ \left\| \sum_{j=0}^{\beta -1}{\frac{\mathbf{1}^{\top}}{n}\mathbf{y}_{\beta k+j}} \right\| ^2 \right] }}
\\
&\quad -2\gamma \beta \underset{=:S_3}{\underbrace{\mathbb{E} \left[ \left< \bar{x}_{\beta k}-x^*,\frac{\mathbf{1}^{\top}}{n}\nabla F_{\beta k} \right> \right] }}.
\end{aligned}
\end{equation}

Next, we bound the terms $S_2$ and $S_3$ in the above equation, respectively.
For $S_2$, by \eqref{Eq_z_update} we have
\begin{equation}\label{Eq_norm_of_tra_error}
\begin{aligned}
&\mathbb{E} \left[ \left\| \sum_{j=0}^{\beta -1}{\frac{\mathbf{1}^{\top}}{n}\mathbf{y}_{\beta k+j}} \right\| ^2 \right] 
\\
&\leqslant \frac{\beta \sigma ^2}{n}+\mathbb{E} \left[ \left\| \sum_{j=0}^{\beta -1}{\frac{1}{n}\sum_{i=1}^n{\nabla f_i\left( x_{i,\beta k} \right)}} \right\| ^2 \right] 
\\
&\leqslant \frac{\beta \sigma ^2}{n}+\frac{2\beta ^{2}L^2}{n}\mathbb{E} \left[ \left\| \tilde{\mathbf{x}}_{\beta k} \right\| ^2 \right] 
+4\beta ^{2}L\mathbb{E} \left[ f\left( \bar{x}_{\beta k} \right) -f\left( x^* \right) \right] ,
\end{aligned}
\end{equation}
where we have used Assumption~\ref{Ass_bounded_var} on bounded stochastic gradient variance and Assumption~\ref{Ass_smooth} on smoothness of the objectives.
For $S_3$, using the convexity and smoothness of $f_i$ in Assumptions~\ref{Ass_convexity} and \ref{Ass_smooth}, we have
\begin{equation}
\begin{aligned}
&\mathbb{E} \left[ \left< \bar{x}_{\beta k}-x^*,\frac{\mathbf{1}^{\top}}{n}\nabla F_{\beta k} \right> \right] 
\\
&=\frac{1}{n}\sum_{i=1}^n{\mathbb{E} \left[ \left< x_{i,\beta k}-x^*,\nabla f_i\left( x_{i,\beta k} \right) \right> \right]}
\\
&\quad-\frac{1}{n}\sum_{i=1}^n{\mathbb{E} \left[ \left< x_{i,\beta k}-\bar{x}_{\beta k},\nabla f_i\left( x_{i,\beta k} \right) \right> \right]}
\\
&\geqslant \mathbb{E} \left[ f\left( \bar{x}_{\beta k} \right) -f\left( x^* \right) \right] 
+\frac{\mu}{2n}\mathbb{E} \left[ \left\| \mathbf{x}_{\beta k}-\mathbf{1}x^* \right\| ^2 \right] 
\\
&\quad-\frac{L}{2n}\mathbb{E} \left[ \left\| \tilde{\mathbf{x}}_{\beta k} \right\| ^2 \right] .
\end{aligned}
\end{equation}

Substituting $S_2$ and $S_3$ in \eqref{Proof_lem1_1}, and noticing that
\begin{equation}
\begin{aligned}
&-\mathbb{E} \left[ \left\| \mathbf{x}_{\beta k}-\mathbf{1}x^* \right\| ^2 \right] 
\\
&=-\mathbb{E} \left[ \left\| \tilde{\mathbf{x}}_{\beta k} \right\| ^2 \right] -n\mathbb{E} \left[ \left\| \bar{x}_{\beta k}-x^* \right\| ^2 \right] 
\\
&\quad-2\sum_{i=1}^n{\left( \mathbb{E} \left[ \left< x_{i,\beta k}-\bar{x}_{\beta k}, \bar{x}_{\beta k}-x^* \right> \right] \right)}
\\
&\leqslant \mathbb{E} \left[ \left\| \tilde{\mathbf{x}}_{\beta k} \right\| ^2 \right] -\frac{n}{2}\mathbb{E} \left[ \left\| \bar{x}_{\beta k}-x^* \right\| ^2 \right] ,
\end{aligned}
\end{equation}
we get
\begin{equation}
\begin{aligned}
&\mathbb{E} \left[ \left\| \bar{x}_{\beta \left( k+1 \right)}-x^* \right\| ^2 \right] 
\\
&\leqslant \left( 1-\frac{\mu \beta \gamma}{2} \right) \mathbb{E} \left[ \left\| \bar{x}_{\beta k}-x^* \right\| ^2 \right] +\frac{\gamma ^2\beta \sigma ^2}{n}
\\
&\quad +\left( \frac{2\gamma ^2\beta ^{2}L^2}{n}+\frac{2\gamma \beta L}{n} \right) \mathbb{E} \left[ \left\| \tilde{\mathbf{x}}_{\beta k} \right\| ^2 \right] 
\\
&\quad-\left( 2\gamma \beta -4\gamma ^2\beta ^{2}L \right) \mathbb{E} \left[ f\left( \bar{x}_{\beta k} \right) -f\left( x^* \right) \right].
\end{aligned}
\end{equation}
Then, letting the stepsize $\gamma \leqslant \frac{1}{4\beta L}$, we complete the proof.
\end{proof}

Next, we present the descent lemma for the gradient norm in the nonconvex case, which plays a key role in establishing the convergence result in \eqref{Eq_convergence_FlexGT_nc}.

\begin{Lem}\label{Lem_descent_FlexGT}
    Suppose Assumptions~\ref{Ass_smooth}--\ref{Ass_graph} hold. Let the stepsize satisfy $\gamma \leqslant \frac{1}{4L\beta }$. Then, we have for all $k \geqslant 0$, 
\begin{equation}\label{Eq_descent_lem}
\begin{aligned}
&\frac{1}{K}\sum_{k=0}^{K-1}{\mathbb{E} \left[ \left\| \nabla f\left( \bar{x}_{\beta k} \right) \right\| ^2 \right]}
\\
&\leqslant \frac{4\left( \mathbb{E} \left[ f\left( \bar{x}_0 \right) \right] -\mathbb{E} \left[ f^* \right] \right)}{\gamma \beta K}
\\
&\quad+\frac{4L^2}{n}\frac{1}{K}\sum_{k=0}^{K-1}{\mathbb{E} \left[ \left\| \tilde{\mathbf{x}}_{\beta k} \right\| ^2 \right]}+\frac{2\gamma L}{n}\sigma ^2.
\end{aligned}
\end{equation}
\end{Lem}

\begin{proof}
By the smoothness of the objective function in Assumption~\ref{Ass_smooth}, we get
\begin{equation}
\begin{aligned}
f\left( \bar{x}_{\beta \left( k+1 \right)} \right) 
&\leqslant f\left( \bar{x}_{\beta k} \right) +\left< \nabla f\left( \bar{x}_{\beta k} \right) ,\bar{x}_{\beta \left( k+1 \right)}-\bar{x}_{\beta k} \right> 
\\
&\quad+\frac{L}{2}\left\| \bar{x}_{\beta \left( k+1 \right)}-\bar{x}_{\beta k} \right\| ^2.
\end{aligned}
\end{equation}
Then, for the inner-product term on the right, by the recursion of decision variable $x$ in \eqref{Eq_x_update} and noticing that $\mathbb{E} \left[ G_{\beta k+j} \right] =\nabla F_{\beta k}$, we have
\begin{equation}
\begin{aligned}
&\mathbb{E} \left[ \left< \nabla f\left( \bar{x}_{\beta k} \right) ,\bar{x}_{\beta \left( k+1 \right)}-\bar{x}_{\beta k} \right> \right] 
\\
&=\mathbb{E} \left[ \left< \nabla f\left( \bar{x}_{\beta k} \right) ,-\gamma \sum_{j=0}^{\beta -1}{\frac{\mathbf{1}^{\top}}{n}\mathbf{y}_{\beta k+j}} \right> \right] 
\\
&\leqslant -\frac{\gamma}{2}\beta \mathbb{E} \left[ \left\| \nabla f\left( \bar{x}_{\beta k} \right) \right\| ^2 \right] 
\\
&\quad+\frac{\gamma}{2}\mathbb{E} \left[ \left\| \sum_{j=0}^{\beta -1}{\frac{\mathbf{1}^{\top}}{n}\left( \nabla F_{\beta k}-\nabla F\left( \mathbf{1}\bar{x}_{\beta k} \right) \right)} \right\| ^2 \right] 
\\
&\leqslant -\frac{\gamma \beta }{2}\mathbb{E} \left[ \left\| \nabla f\left( \bar{x}_{\beta k} \right) \right\| ^2 \right] +\frac{\gamma \beta L^2}{2n}\mathbb{E} \left[ \left\| \tilde{\mathbf{x}}_{\beta k} \right\| ^2 \right] ,
\end{aligned}
\end{equation}

Then, based on the upper bound of the GT error obtained in \eqref{Eq_norm_of_tra_error}, we have
\begin{equation}
\begin{aligned}
&\mathbb{E} \left[ f\left( \bar{x}_{\beta \left( k+1 \right)} \right) \right] 
\\
&\leqslant \mathbb{E} \left[ f\left( \bar{x}_{\beta k} \right) \right] -\left( \frac{\gamma \beta }{2}-\gamma ^2\beta ^{2}L \right) \mathbb{E} \left[ \left\| \nabla f\left( \bar{x}_{\beta k} \right) \right\| ^2 \right] 
\\
&\quad +\left( \frac{\gamma \beta L^2}{2n}+\frac{\gamma ^2\beta ^{2}L^3}{n} \right) \mathbb{E} \left[ \left\| \tilde{\mathbf{x}}_{\beta k} \right\| ^2 \right] +\frac{\gamma ^2\beta L}{2n}\sigma ^2.
\end{aligned}
\end{equation}
Letting $\gamma \leqslant \frac{1}{8L\beta }$, we get
\begin{equation}
\begin{aligned}
\mathbb{E} \left[ \left\| \nabla f\left( \bar{x}_{\beta k} \right) \right\| ^2 \right] 
&\leqslant \frac{4\mathbb{E} \left[ f\left( \bar{x}_{\beta k} \right) -f\left( \bar{x}_{\beta \left( k+1 \right)} \right) \right]}{\gamma \beta }
\\
&\quad+\frac{4L^2}{n}\mathbb{E} \left[ \left\| \tilde{\mathbf{x}}_{\beta k} \right\| ^2 \right] +\frac{2\gamma L}{n}\sigma ^2.
\end{aligned}
\end{equation}
    Summing the above inequality from $0$ to $k-1$, we obtained the targeted results.
\end{proof}

\subsection{Proof of Theorem~\ref{Thm_FlexGT}}\label{Sec_proof_Thm_FlexGT}

With Lemmas~\ref{Lem_opt_gap} and~\ref{Lem_descent_FlexGT} for optimality gap and objective function decrease, Lemmas~\ref{Lem_cons_err} and~\ref{Lem_cons_err_sum_FlexGT} for consensus errors, and Lemmas~\ref{Lem_tra_err} and~\ref{Lem_trac_err_sum_FlexGT} for tracking errors, we are now ready to prove Theorem~\ref{Thm_FlexGT}.

\textbf{Strongly convex case.} Recall the Lyapunov function \eqref{Def_Lyapunov_func} with the coefficients set as follows:
\begin{equation}\label{Eq_c1_c2}
\begin{aligned}
c_x=\frac{16\gamma \beta L}{n\left( 1-\bar{\rho}_W \right)},\quad c_y=\frac{256\gamma ^3\beta ^{3}L\bar{\rho}_W}{n\left( 1-\bar{\rho}_W \right) ^3}.
\end{aligned}
\end{equation}
Then, by substituting the results obtained from Lemmas~\ref{Lem_opt_gap},~\ref{Lem_cons_err} and~\ref{Lem_tra_err} into the Lyapunov function \eqref{Def_Lyapunov_func}, and noticing that
\[\left\| \nabla f\left( \bar{x}_{\beta k} \right) \right\| ^2\leqslant 2L\left( f\left( \bar{x}_{\beta k} \right) -f\left( x^* \right) \right), \]
we can get the contraction property for the Lyapunov function as follows:
\begin{equation}
\begin{aligned}
\mathbb{E} \left[ V_{\beta \left( k+1 \right)} \right] &\leqslant \left( 1-\max \left\{ \frac{\mu \beta \gamma}{2}, \frac{1-\bar{\rho}_W}{8} \right\} \right) \mathbb{E} \left[ V_{\beta k} \right] 
\\
&\quad +e_1\mathbb{E} \left[ f\left( \bar{x}_{\beta k} \right) -f\left( x^* \right) \right] 
\\
&\quad +e_2\mathbb{E} \left[ \left\| \tilde{\mathbf{x}}_{\beta k} \right\| ^2 \right] +e_3\mathbb{E} \left[ \left\| \tilde{\mathbf{y}}_{\beta k} \right\| ^2 \right] 
\\
&\quad +\frac{\gamma ^2\beta \sigma ^2}{n}+c_x\frac{8n\gamma ^2\beta ^{2}\bar{\rho}_W}{1-\bar{\rho}_W}\sigma ^2+c_y6n\sigma ^2,
\end{aligned}
\end{equation}
where
\begin{equation*}
\begin{aligned}
e_1&:=\frac{192n\gamma ^2\beta ^{2}L^3\bar{\rho}_W}{1-\bar{\rho}_W}c_y-\gamma \beta ,
\\
e_2&:=\frac{3\gamma \beta L}{n}+\frac{18\bar{\rho}_WL^2}{1-\bar{\rho}_W}c_y-\frac{3\left( 1-\bar{\rho}_W \right)}{8}c_x,
\\
e_3&:=\frac{4\gamma ^2\beta ^{2}\bar{\rho}_W}{1-\bar{\rho}_W}c_x-\frac{1-\bar{\rho}_W}{4}c_y.
\end{aligned}
\end{equation*}
Letting $e_1\leqslant -\gamma \beta /2$ and $e_2, e_3\leqslant 0$ with the stepsize $\gamma$ satisfying 
\begin{equation}\label{Eq_stepsize_sc}
\begin{aligned}
\gamma \leqslant \min \left\{ \frac{1}{4\sqrt{2}\beta L}, \frac{1-\bar{\rho}_W}{18\beta L\sqrt{\bar{\rho}_W}}, \frac{\left( 1-\bar{\rho}_W \right) ^2}{40\beta L\bar{\rho}_W} \right\} ,
\end{aligned}
\end{equation}
we obtain the convergence result for strongly convex case as shown in \eqref{Eq_convergence_FlexGT_sc}.

\textbf{Convex case.}  Under Assumption~\ref{Ass_convexity} with $\mu=0$ and  $x^*$ as one of the optimal solutions. Noticing that $e_1\leqslant -\gamma \beta /2$, we can rewrite the contraction property of the optimality gap in Lemma~\ref{Lem_opt_gap} as follows:
\begin{equation}
\begin{aligned}
&\mathbb{E} \left[ \left\| \bar{x}_{\beta \left( k+1 \right)}-x^* \right\| ^2 \right] 
\\
&\leqslant \mathbb{E} \left[ \left\| \bar{x}_{\beta k}-x^* \right\| ^2 \right] +\frac{\gamma ^2\beta \sigma ^2}{n}
\\
&\quad +\frac{3\gamma \beta L}{n}\mathbb{E} \left[ \left\| \tilde{\mathbf{x}}_{\beta k} \right\| ^2 \right] -\gamma \beta \mathbb{E} \left[ f\left( \bar{x}_{\beta k} \right) -f^*\left( x^* \right) \right] .
\end{aligned} 
\end{equation}
Following the approach used in the strongly convex case, we apply the contracted Lyapunov function with the same values of $c_x$ and $c_y$ as in \eqref{Eq_c1_c2}.
By applying Lemmas~\ref{Lem_cons_err} and~\ref{Lem_tra_err} and ensuring that the stepsize satisfy \eqref{Eq_stepsize_sc},
we have
\begin{equation}
\begin{aligned}
\mathbb{E} \left[ f\left( \bar{x}_{\beta k} \right) -f\left( x^* \right) \right] 
&\leqslant \frac{2\left( \mathbb{E} \left[ V_{\beta k} \right] -\mathbb{E} \left[ V_{\beta \left( k+1 \right)} \right] \right)}{\gamma \beta }
\\
&\quad +\frac{2\gamma \sigma ^2}{n}+\frac{3328\gamma ^2\beta ^{2}L\bar{\rho}_W}{\left( 1-\bar{\rho}_W \right) ^3}\sigma ^2.
\end{aligned}
\end{equation}
Summing the above inequality from $k = 0$ to $K - 1$ and rearranging the terms, we obtain the result for the convex case as given in \eqref{Eq_convergence_FlexGT_gc}.

\textbf{Nonconvex case.} 
By combining the obtained Lemma~\ref{Lem_cons_err_sum_FlexGT} and Lemma~\ref{Lem_trac_err_sum_FlexGT} regarding the accumulated consensus error and GT error, respectively, we can further obtain
\begin{equation}
\begin{aligned}
&\left( 1-\frac{768\gamma ^2\beta ^{2}L^2\bar{\rho}_{W}^{2}}{\left( 1-\bar{\rho}_W \right) ^4} \right) \frac{1}{K}\sum_{k=0}^{K-1}{\mathbb{E} \left[ \left\| \tilde{\mathbf{x}}_{\beta k} \right\| ^2 \right]}
\\
&\leqslant \frac{2\mathbb{E} \left[ \left\| \tilde{\mathbf{x}}_0 \right\| ^2 \right]}{\left( 1-\bar{\rho}_W \right) K}+\frac{32\gamma ^2\beta ^{2}\bar{\rho}_W\mathbb{E} \left[ \left\| \tilde{\mathbf{y}}_0 \right\| ^2 \right]}{\left( 1-\bar{\rho}_W \right) ^3K}
\\
&\quad +\frac{16n\gamma ^2\beta ^{2}\bar{\rho}_W}{\left( 1-\bar{\rho}_W \right) ^2}\sigma ^2+\frac{192n\gamma ^2\beta ^{2}\bar{\rho}_W}{\left( 1-\bar{\rho}_W \right) ^3}\sigma ^2
\\
&\quad +\frac{3072n\gamma ^4\beta ^{4}L^2\bar{\rho}_{W}^{2}}{\left( 1-\bar{\rho}_W \right) ^4}\frac{1}{K}\sum_{k=0}^{K-1}{\mathbb{E} \left[ \left\| \nabla f\left( \bar{x}_{\beta k} \right) \right\| ^2 \right]}.
\end{aligned}
\end{equation}
Letting the stepsize satisfy 
\begin{equation}\label{Eq_stepsize_accumulated}
\gamma \leqslant \min \left\{ \frac{1}{4\beta L}, \frac{1-\bar{\rho}_W}{14\beta L\sqrt{\bar{\rho}_W}}, \frac{\left( 1-\bar{\rho}_W \right) ^2}{40\beta L\bar{\rho}_W} \right\},
\end{equation}
we get for all $K \geqslant 1$,
\begin{equation}\label{Eq_con_err_decoupled}
\begin{aligned}
\frac{1}{K}\sum_{k=0}^{K-1}{\mathbb{E} \left[ \left\| \tilde{\mathbf{x}}_{\beta k} \right\| ^2 \right]}
&\leqslant \frac{C_{0}}{\left( 1-\bar{\rho}_W \right) K}+\frac{416n\gamma ^2\beta ^{2}\bar{\rho}_W}{\left( 1-\bar{\rho}_W \right) ^3}\sigma ^2
\\
&\quad +4n\gamma ^2\beta ^{2}\frac{1}{K}\sum_{k=0}^{K-1}{\mathbb{E} \left[ \left\| \nabla f\left( \bar{x}_{\beta k} \right) \right\| ^2 \right]},
\end{aligned}
\end{equation}
where
\begin{equation}\label{Eq_C_0_FlexGT}
C_0:=4\mathbb{E} \left[ \left\| \tilde{\mathbf{x}}_0 \right\| ^2 \right] +\frac{64\gamma ^2\beta ^{2}\bar{\rho}_W\mathbb{E} \left[ \left\| \tilde{\mathbf{y}}_0 \right\| ^2 \right]}{\left( 1-\bar{\rho}_W \right) ^2}.
\end{equation}
Then, by incorporating \eqref{Eq_con_err_decoupled} into \eqref{Eq_descent_lem} and 
setting the initial consensus error in \eqref{Eq_C_0_FlexGT} to zero, i.e., $x_{i,0}=x_{j,0}, \forall i,j\in[n]$, we complete the proof.

\subsection{Proof of Corollary \ref{Col_complexity_FlexGT}}
\label{Proof_Cor_1}
Denote by $\hat{\gamma}$ the obtained upper bound of the stepsize in \eqref{Eq_stepsize_global} and 
\begin{equation*}
\varUpsilon _1:=\frac{\beta \sigma ^2}{n\upsilon  _0},\quad \varUpsilon _2:=\frac{\beta ^{3}L\bar{\rho}_W\sigma ^2}{\upsilon  _0\left( 1-\bar{\rho}_W \right) ^3},
\end{equation*}
the constant part of the centralized error and the topology effect part in \eqref{Eq_convergence_FlexGT_sc}--\eqref{Eq_convergence_FlexGT_nc}, respectively, both of which are controlled by the stepsize $\gamma$. Here,
$\upsilon _0$ denotes either $\mathbb{E} \left[ V_0 \right] $ or $\mathbb{E} \left[ f\left( \bar{x}_0 \right) -f^* \right]$, depending on the specific case under consideration.
Then, for the strongly convex case, tune the stepsize between
\begin{equation*}
\left\{ \frac{2\ln \left( \max \left\{ 2,\beta ^{2}\mu ^2\hat{K}/\varUpsilon _1, \beta ^{3}\mu ^3\hat{K}^2/\varUpsilon _2 \right\} \right)}{\beta \mu \hat{K}}, \,\,\hat{\gamma} \right\}
\end{equation*}
to ensure that the centralized error and topology effect part match the linear part in terms of the number of iterations $\hat{K}$. Specifically, if $\hat{\gamma}$ is smaller, 
we have
\begin{equation}\label{Eq_proof_comp}
\begin{aligned}
\mathbb{E} \left[ V_{\beta k} \right] &\leqslant \exp \left( -\frac{\mu \left( 1-\bar{\rho}_W \right) ^2}{L}K \right) \mathbb{E} \left[ V_0 \right] 
\\
&\quad+\left( \frac{\varUpsilon _1}{\beta ^2\mu ^2\hat{K}}+\frac{\varUpsilon _2}{\beta ^3\mu ^3\hat{K}^2} \right) \mathbb{E} \left[ V_0 \right] 
\end{aligned}
\end{equation}
Otherwise,
noticing that
\begin{equation}
\begin{aligned}
&\exp \left( -\min \left\{ \frac{\mu \beta \gamma}{2},\frac{1-\bar{\rho}_W}{8} \right\} K \right) \mathbb{E} \left[ V_0 \right] 
\\
&\leqslant \exp \left( -\ln \left( \max \left\{ 2,\frac{\beta ^2\mu ^2\hat{K}}{\varUpsilon _1},\frac{\beta ^3\mu ^3\hat{K}^2}{\varUpsilon _2} \right\} \right) \right) \mathbb{E} \left[ V_0 \right] 
\\
&\leqslant \left( \frac{\varUpsilon _1}{\beta ^2\mu ^2\hat{K}}+\frac{\varUpsilon _2}{\beta ^3\mu ^3\hat{K}^2} \right) \mathbb{E} \left[ V_0 \right],
\end{aligned}
\end{equation}
we get the same result as \eqref{Eq_proof_comp}. Letting the right-hand side of \eqref{Eq_proof_comp} further smaller than or equal to $\epsilon$, we obtain the complexity result in \eqref{Eq_complexity_FlexGT_sc}.

Similarly, for the convex and nonconvex cases, tuning the stepsize between
\begin{equation*}
\left\{ \frac{1}{\sqrt{\varUpsilon _1\hat{K}}},\,\, \frac{1}{\sqrt[3]{\varUpsilon _2\hat{K}}},\,\, \hat{\gamma} \right\},
\end{equation*}
we obtain the complexity results in \eqref{Eq_complexity_FlexGT_gc} and \eqref{Eq_complexity_FlexGT_nc}. 

\subsection{Proof of Corollary~\ref{Cor_Acc-FlexGT}}
\label{Proof_cor_2}
The proof is similar to that of Corollary~\ref{Col_complexity_FlexGT}. By Lemma~\ref{Lem_communication_protocol}, we have $1-\bar{\rho}_W \geqslant 1/2$, and the topology effect part $\varUpsilon _2$ can be dominated when $\bar{\rho}_W$ is sufficiently small. Specifically, let $\hat{\gamma}$ denote the obtained upper bound of the stepsize in \eqref{Eq_stepsize_global}, and tune the stepsize $\gamma$ between
\begin{equation*}
\left\{ \frac{2\ln \left( \max \left\{ 2,\beta ^{2}\mu ^2\hat{K}/\varUpsilon _1 \right\} \right)}{\beta \mu \hat{K}}, \,\, \hat{\gamma} \right\} 
\end{equation*}
for the strongly convex case; and
\begin{equation*}
\left\{ \frac{1}{\sqrt{\varUpsilon _1\hat{K}}}, \,\, \hat{\gamma} \right\}
\end{equation*}
for the convex and nonconvex cases to ensure that the centralized error part $\varUpsilon _1$ matches the linear or sublinear term. Then, we obtain the target complexity results in \eqref{Eq_complexity_Acc_FlexGT_sc}, \eqref{Eq_complexity_Acc_FlexGT_gc}, and \eqref{Eq_complexity_Acc_FlexGT_nc}.

\section{Numerical Experiments}\label{Sec_experiments}
In this section, we report several numerical experiments to validate the theoretical findings of FlexGT and Acc-FlexGT on solving the distributed stochastic optimization problem \eqref{Prob}, using both synthetic and real-world datasets. We compare FlexGT and Acc-FlexGT with several baseline methods, including DSGD \cite{lian2017can}, DSGT \cite{pu2021distributed}, DFL \cite{liu2022decentralized}, K-GT \cite{liu2025decentralized}, and SS-DSGT \cite{didouble}.

\subsection{Synthetic Example}

We consider the following distributed ridge regression problem over a network of $n=20$ nodes:
\begin{equation}\label{Prob_RR}
\begin{aligned}
\underset{x\in \mathbb{R}^p}{\min}f\left( x \right) =\frac{1}{n}\sum_{i=1}^n{\underset{=:f_i}{\underbrace{\left( \mathbb{E} _{v_i}\left[ \left( h_{i}^{T}x-\nu _i \right) ^2+\frac{\mu}{2}\left\| x \right\| ^2 \right] \right) }}},
\end{aligned}
\end{equation}
where $\mu>0$ is the regularization parameter, $h_i\in \left[ 0,1 \right] ^p$ denotes the feature parameters of node $i$ with dimension $p=10$, and $v_i\sim \mathcal{N} \left( \bar{v}_i,\sigma ^2 \right)$ with $\bar{v}_i\in \left[ 0,1 \right]$. Therefore, the algorithm can obtain an unbiased noisy gradient $\nabla f_i\left( x_{i,t} \right) +\delta _{i,t}$ with $\delta _{i,t}\sim \mathcal{N} \left( 0,\sigma^2 \right) $ at each iteration $t$.
Moreover, we set the stepsize according to the choice of $\alpha $ and $\beta $ as suggested by Theorem~\ref{Thm_FlexGT}, i.e., $\gamma =c\left( 1-\bar{\rho} _{W} \right) ^2/\left(\bar{\rho} _{W} \beta L \right) $, where $c=10$ is a constant and the Lipschitz constant is set to $L=1$. We note that there is heterogeneity among the objectives $f_i$ of nodes due to the differences in $\{ h_i \} _{i=1}^{n}$.

\begin{figure}[!tb]
    \centering
    \subfloat[$\mu=0.001$]
    {
        \begin{minipage}[t]{0.24\textwidth}
            \centering
            \includegraphics[width=\textwidth]{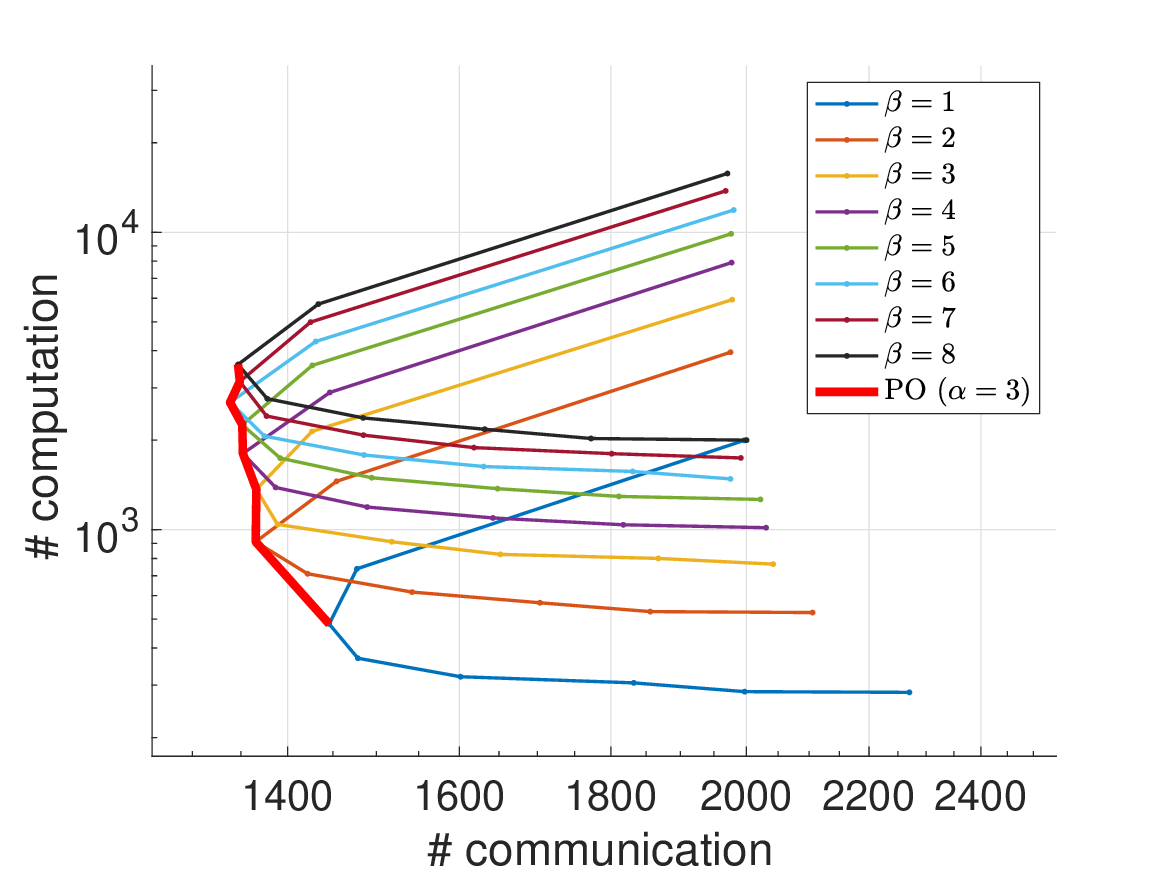}
        \end{minipage}%
    }
    \subfloat[$\mu=0.0001$]
    {
        \begin{minipage}[t]{0.24\textwidth}
            \centering
            \includegraphics[width=\textwidth]{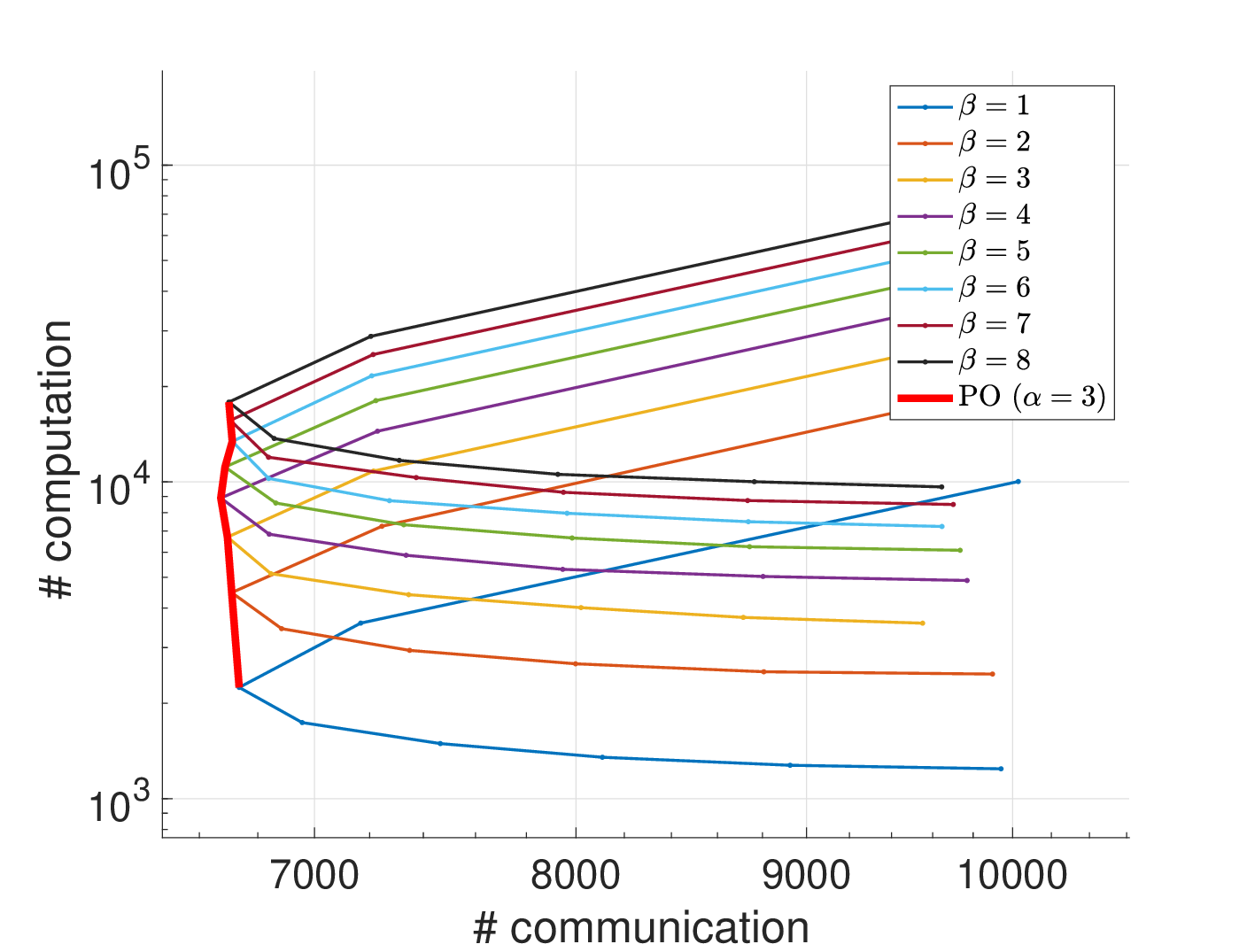}
        \end{minipage}
    }   
    \caption{Communication and computational complexity of FlexGT algorithm to achieve $\epsilon=10^{-3}$ accuracy with $\alpha=1,2,\dots,8$ and $\beta=1,2,\dots,8$ on synthetic data. Each node is in an exponential graph of $n=20$ nodes with 5 neighbors.}
    \label{Fig_exp_pareto}
\end{figure}

In Fig.~\ref{Fig_exp_pareto}, we present the communication and computational complexity of FlexGT to achieve an accuracy of $\epsilon = 10^{-3}$ with $\alpha = 1,2,\dots,8$ and $\beta = 1,2,\dots,8$. The experiment was run 100 times for each parameter setting to obtain an average complexity. The results demonstrate that, for a fixed number of communication steps $\alpha$, increasing the number of computation steps per round $\beta$ reduces the communication complexity but increases the computational complexity (moving towards the upper left). However, when $\beta$ is fixed and $\alpha$ is increased, the computational complexity decreases initially, while the communication complexity first decreases and then increases.
It is worth noting that the Pareto-optimal solutions (bold red lines) are achieved according to the value of $\alpha$ suggested by Corollary~\ref{Cor_Acc-FlexGT}. 
For smaller $\mu$, more communication and computation are required, as shown in the plot on the right. These observations validate our theoretical findings and illustrate the trade-off between communication and computation in distributed optimization algorithms.

\begin{figure}[!tb]
    \centering
    \subfloat[$\left| \mathcal{N} _i \right|=5, \rho _W=0.74$]
    {
        \begin{minipage}[t]{0.24\textwidth}
            \centering
            \includegraphics[width=\textwidth]{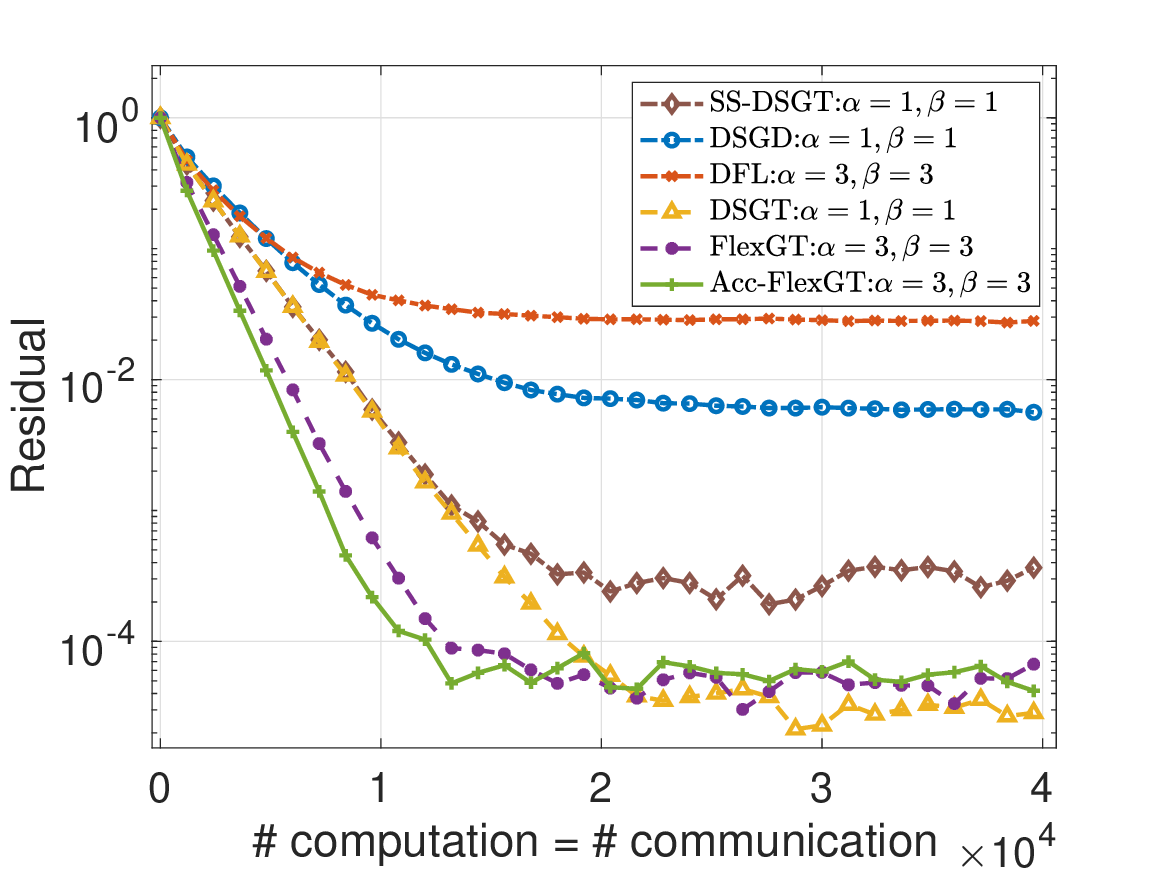}
        \end{minipage}%
    }
    \subfloat[$\left| \mathcal{N} _i \right|=3, \rho _W=0.88$]
    {
        \begin{minipage}[t]{0.24\textwidth}
            \centering
            \includegraphics[width=\textwidth]{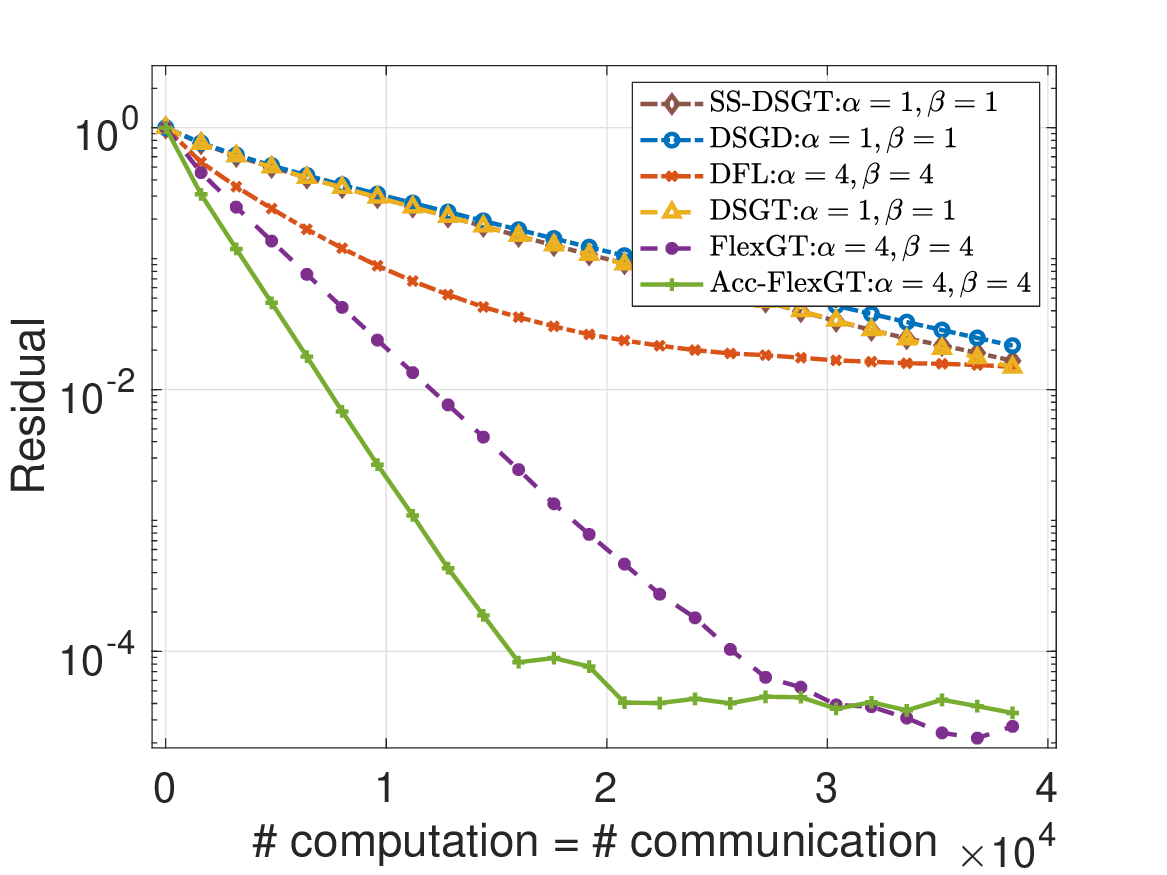}
        \end{minipage}
    }   
    \caption{Comparison of convergence performance on synthetic data. Each node in an exponential graph of $n=20$ nodes with $\left| \mathcal{N} _i \right|$ neighbors.}
    \label{Fig_exp_RR_convergence}
\end{figure}

In Fig.~\ref{Fig_exp_RR_convergence}, we compare the convergence of the residual $\left| x_k - x^* \right|^2 / \left| x_0 - x^* \right|^2$ for the problem \eqref{Prob_RR} between FlexGT, Acc-FlexGT, and baseline methods. Following the recommendation of Acc-FlexGT, we set the number of communication and computation steps according to \eqref{Eq_d1_Acc_FlexGT_sc} with $\alpha = \beta$.
The results show that both FlexGT and Acc-FlexGT achieve superior convergence rates compared to other algorithms. This advantage is particularly evident in the scenario with weaker topological connectivity, i.e., when $\rho_W$ is larger, as shown in the plot on the right. 

\subsection{MNIST Dataset} 
\begin{figure*}[!tb]
    \centering
    \subfloat
    {
        \begin{minipage}[b]{0.3\textwidth}
            \centering  
            \includegraphics[width=\textwidth]{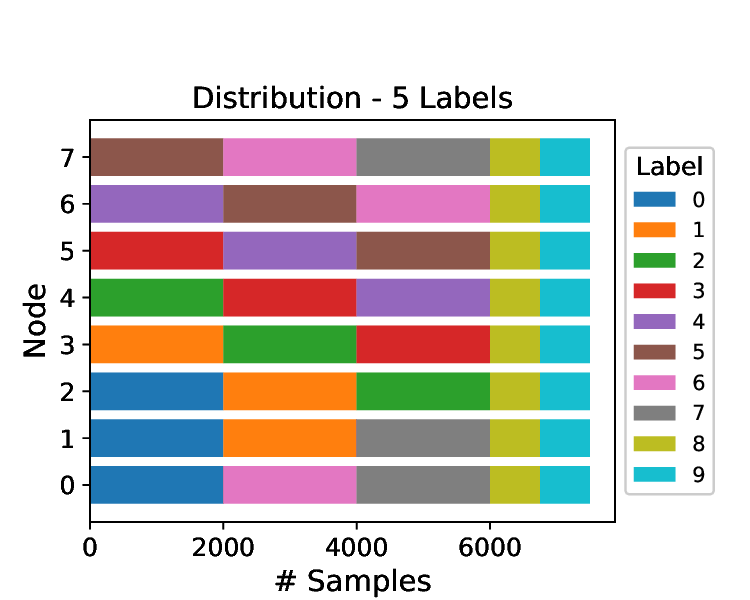}
        \end{minipage}
    }
    \subfloat
    {
        \begin{minipage}[b]{0.3\textwidth}
            \centering  
            \includegraphics[width=\textwidth]{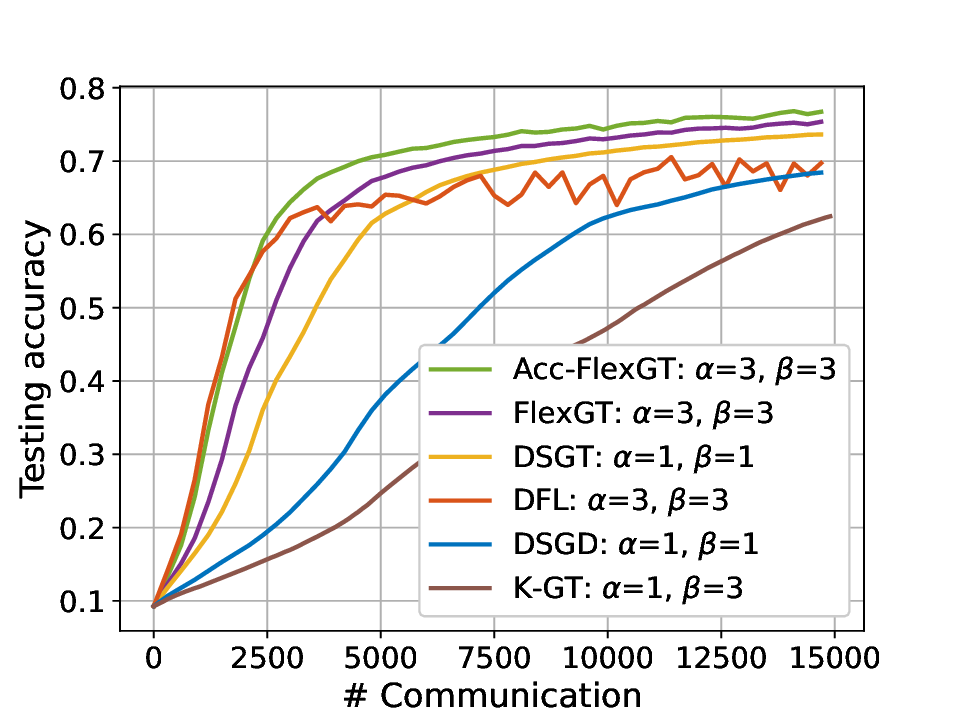}
        \end{minipage}
    }
    \subfloat
    {
        \begin{minipage}[b]{0.3\textwidth}
            \centering  
            \includegraphics[width=\textwidth]{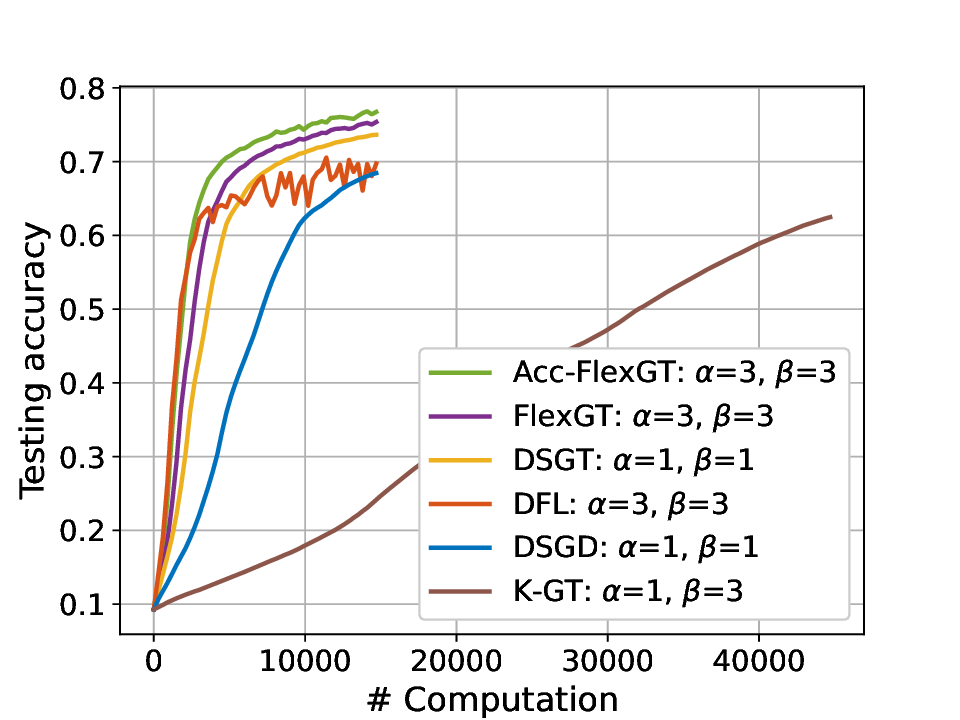}
        \end{minipage}
    }

    \subfloat
    {
        \begin{minipage}[b]{0.3\textwidth}
            \centering  
            \includegraphics[width=\textwidth]{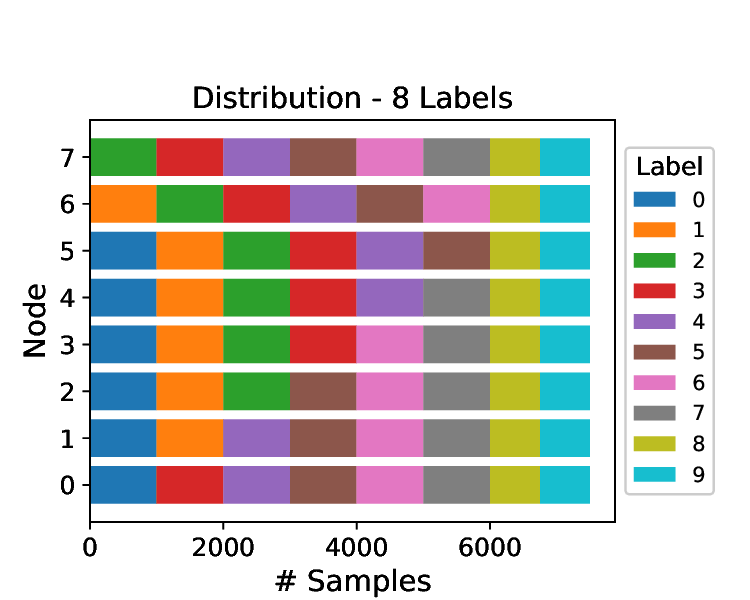}
        \end{minipage}
    }
    \subfloat
    {
        \begin{minipage}[b]{0.3\textwidth}
            \centering  
            \includegraphics[width=\textwidth]{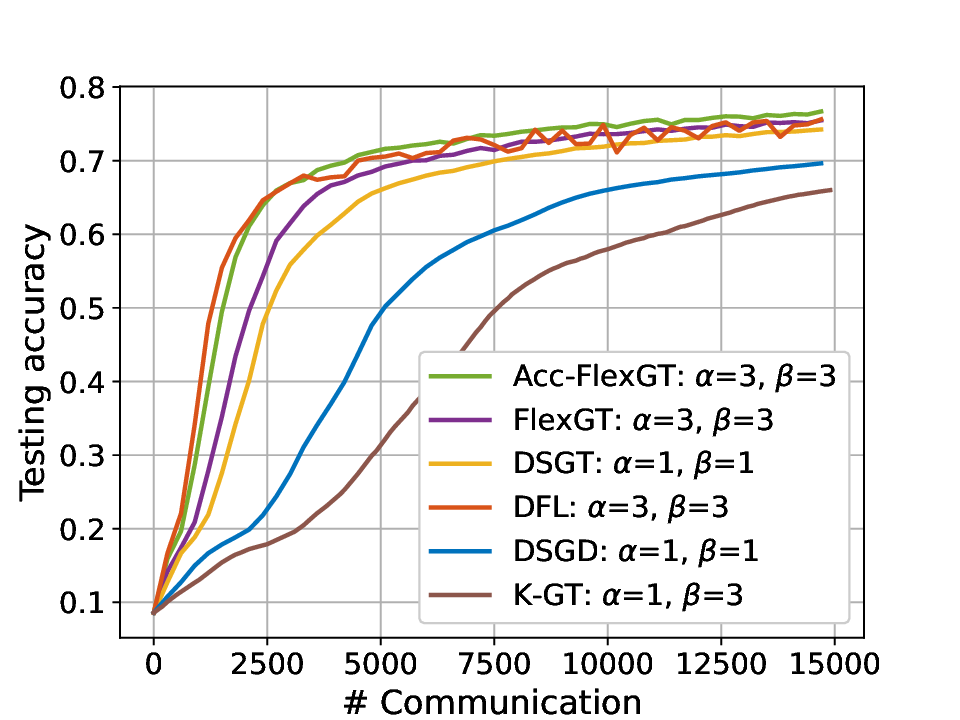}
        \end{minipage}
    }
    \subfloat
    {
        \begin{minipage}[b]{0.3\textwidth}
            \centering  
            \includegraphics[width=\textwidth]{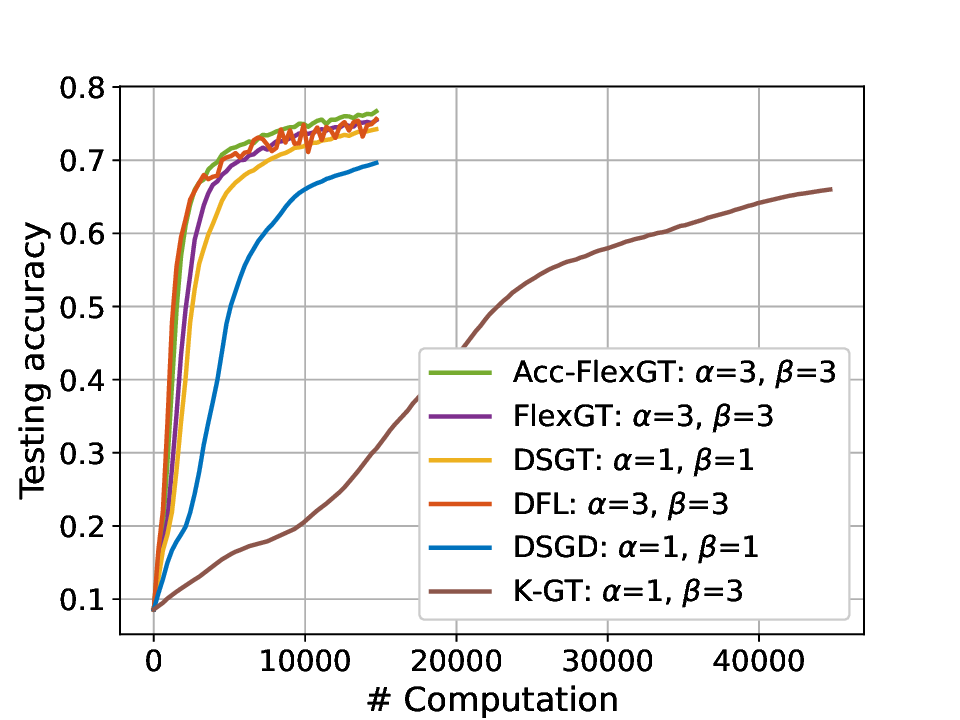}
        \end{minipage}
    }   
    \caption{Comparison of test accuracy against computation and communication steps on the MNIST dataset with heterogeneous label distributions. The top plots represent each node having data from 5 classes (higher heterogeneity), while the bottom plots correspond to each node having data from 8 classes (lower heterogeneity).}
    \label{Fig_training}
\end{figure*}

We conduct distributed training of three-layer neural networks with smooth activation functions for 10-class classification on the MNIST dataset \cite{deng2012mnist}. The original training set is balanced, containing an equal number of samples for each label from 0 to 9. To introduce heterogeneous label distributions, we partition the training data across 8 GPUs connected over a ring graph, where each GPU receives samples corresponding to only 8 or 5 classes, as illustrated in the left column of Fig.~\ref{Fig_training}.
The average Kullback-Leibler divergence \cite{kullback1951information} for these two label distributions is 0.23 and 0.78, respectively, where a higher value indicates greater heterogeneity between each node's local dataset and the overall training set.

Specifically, we evaluate the trained model on the full test dataset and compare the testing accuracy of FlexGT and Acc-FlexGT against baseline methods, including DSGD, DSGT, K-GT, and DFL. As shown in Fig.~\ref{Fig_training}, with $\alpha = \beta = 3$, when the data heterogeneity is high (top plots), FlexGT and Acc-FlexGT outperform the baseline methods in both communication and computational efficiency, with Acc-FlexGT achieving further performance improvements. Notably, the computational complexity of K-GT is worse due to simply skipping communications. These results illustrate the effectiveness of the proposed algorithms.
Additionally, when data heterogeneity is relatively low (bottom plots), the gap between DSGD-based and GT-based algorithms becomes smaller, highlighting the robustness of GT in handling data heterogeneity.

\section{Conclusions}\label{Sec_conclusion}
In this paper, we have proposed FlexGT, a flexible snapshot GT method, along with its accelerated variant, Acc-FlexGT, for solving distributed stochastic optimization problems in non-i.i.d. scenarios. Leveraging a unified convergence analysis framework for strongly convex, convex, and nonconvex objective functions, we showed that FlexGT and Acc-FlexGT achieve linear or sublinear convergence rates with the best-known communication and computational complexity in most settings. In particular, Acc-FlexGT attains optimal iteration complexity and optimal communication complexity in the nonconvex case (up to a logarithmic factor) and improves existing results in the strongly convex case. These findings highlight the inherent trade-off between communication and computation, with Acc-FlexGT achieving a Pareto-optimal solution when connectivity of the graph is available. 

Future work includes exploring single-loop algorithms with randomized communication and computation strategies and evaluating the proposed methods in practical distributed optimization tasks, such as distributed training of large models and multi-robot system applications.

\appendix

This section presents several supporting lemmas for the proof of Theorem~\ref{Thm_FlexGT}. We begin by bounding the client divergence of the decision parameter between two communications.

\begin{Lem}[\textbf{Client divergence}]\label{Lem_tech_client_diver}
Suppose Assumptions~\ref{Ass_smooth}--\ref{Ass_graph} hold. Let the stepsize satisfy $\gamma \leqslant \frac{1}{8\beta L}$.
Then, we have that for all $ k \geqslant 0$
\begin{equation}\label{Eq_Lem_divergence_1}
\begin{aligned}
&\mathbb{E} \left[ \left\| \mathbf{x}_{\beta \left( k+1 \right)}-\mathbf{1}\bar{x}_{\beta k} \right\| ^2 \right] 
\\
&\leqslant 3\mathbb{E} \left[ \left\| \tilde{\mathbf{x}}_{\beta k} \right\| ^2 \right] +16n\gamma ^2\beta ^{2}\sigma ^2
\\
&\quad +8\gamma ^2\beta ^{2}\mathbb{E} \left[ \left\| \tilde{\mathbf{y}}_{\beta k} \right\| ^2 \right] +16n\gamma ^2\beta ^{2}\mathbb{E} \left[ \left\| \nabla f\left( \bar{x}_{\beta k} \right) \right\| ^2 \right] .
\end{aligned}
\end{equation}
\end{Lem}

\begin{proof}
By the $x$-update of FlexGT \eqref{Eq_x_update}, we have
\begin{equation}\label{Eq_lemma9_1}
\begin{aligned}
&\mathbb{E} \left[ \left\| \mathbf{x}_{\beta \left( k+1 \right)}-\mathbf{1}\bar{x}_{\beta k} \right\| ^2 \right] 
\\
&\leqslant \mathbb{E} \left[ \left\| \left( \bar{W}-\mathbf{J} \right) \tilde{\mathbf{x}}_{\beta k}-\gamma \bar{W}\sum_{j=0}^{\beta -1}{\mathbf{y}_{\beta k+j}} \right\| ^2 \right] 
\\
&\leqslant 2\bar{\rho}_W\mathbb{E} \left[ \left\| \tilde{\mathbf{x}}_{\beta k} \right\| ^2 \right] +2\gamma ^2\mathbb{E} \left[ \left\| \sum_{j=0}^{\beta -1}{\mathbf{y}_{\beta k+j}} \right\| ^2 \right] .
\end{aligned}
\end{equation}
Noticing that
\begin{equation}
\begin{aligned}
\mathbb{E} \left[ \left\| \mathbf{y}_{\beta k+t} \right\| ^2 \right] 
&\leqslant 2\mathbb{E} \left[ \left\| \mathbf{y}_{\beta k} \right\| ^2 \right] +4n\sigma ^2
\\
&\leqslant 4\mathbb{E} \left[ \left\| \tilde{\mathbf{y}}_{\beta k} \right\| ^2 \right] +4\mathbb{E} \left[ \left\| \mathbf{1}\bar{y}_{\beta k} \right\| ^2 \right] +4n\sigma ^2
\\
&\leqslant 4\mathbb{E} \left[ \left\| \tilde{\mathbf{y}}_{\beta k} \right\| ^2 \right] +8n\mathbb{E} \left[ \left\| \nabla f\left( \bar{x}_{\beta k} \right) \right\| ^2 \right] 
\\
&\quad +8L^2\mathbb{E} \left[ \left\| \tilde{\mathbf{x}}_{\beta k} \right\| ^2 \right] +4\sigma ^2+4n\sigma ^2.
\end{aligned}
\end{equation}
Plugging the inequality above into \eqref{Eq_lemma9_1} and letting the stepsize satisfy $\gamma \leqslant \frac{1}{4\beta L}$, we complete the proof.
\end{proof}

Then, we bound the consensus error of the decision variables and prove that it contracts at each round, as stated in the following lemma.

\begin{Lem}[\textbf{Consensus error}]\label{Lem_cons_err}
Suppose Assumptions~\ref{Ass_smooth}--\ref{Ass_graph} hold. Let the stepsize satisfy $\gamma \leqslant \frac{1}{4\beta L}$. Then, we have for all $ k \geqslant 0$, 
\begin{equation}
\begin{aligned}
\mathbb{E} \left[ \left\| \tilde{\mathbf{x}}_{\beta \left( k+1 \right)} \right\| ^2 \right] 
&\leqslant \frac{1+\bar{\rho}_W}{2}\mathbb{E} \left[ \left\| \tilde{\mathbf{x}}_{\beta k} \right\| ^2 \right] 
\\
&\quad+\frac{4\gamma ^2\beta ^{2}\bar{\rho}_W}{1-\bar{\rho}_W}\mathbb{E} \left[ \left\| \tilde{\mathbf{y}}_{\beta k} \right\| ^2 \right] +\frac{8n\gamma ^2\beta ^{2}\bar{\rho}_W}{1-\bar{\rho}_W}\sigma ^2.
\end{aligned}
\end{equation}
\end{Lem}

\begin{proof}
According to the update rule of the decision variable in \eqref{Eq_x_update}, we have
\begin{equation}
\begin{aligned}
\mathbb{E} \left[ \left\| \tilde{\mathbf{x}}_{\beta \left( k+1 \right)} \right\| ^2 \right] 
&\leqslant \frac{1+\bar{\rho} _{W}}{2}\mathbb{E} \left[ \left\| \tilde{\mathbf{x}}_{\beta k} \right\| ^2 \right] 
\\
&\quad +\frac{\gamma ^2\beta \left( 1+\bar{\rho} _{W} \right) \bar{\rho} _{W}}{1-\bar{\rho} _{W}}\sum_{t=0}^{\beta -1}{\mathbb{E} \left[ \left\| \tilde{\mathbf{y}}_{\beta k+t}\right\| ^2 \right]},
\end{aligned}
\end{equation}
where we have used Young's inequality  \cite{young1912classes}. 
Noticing that 
\begin{equation}\label{Eq_y_inner_loop}
\mathbb{E} \left[ \left\| \tilde{\mathbf{y}}_{\beta k+t} \right\| ^2 \right] \leqslant 2\mathbb{E} \left[ \left\| \tilde{\mathbf{y}}_{\beta k} \right\| ^2 \right] +4n\sigma ^2,
\end{equation}
and letting the stepsize $\gamma \leqslant \frac{1}{4\beta L}$, we complete the proof.
\end{proof}

Based on Lemma~\ref{Lem_tech_client_diver}, we can further bound the gradient tracking error with a similar contraction property.

\begin{Lem}[\textbf{GT error}]\label{Lem_tra_err}
Suppose Assumptions~\ref{Ass_smooth}--\ref{Ass_graph} hold. Let the stepsize satisfy
\begin{equation}\label{Eq_step_tracking_err}
\gamma \leqslant \min \left\{ \frac{1}{4\beta L}, \frac{1-\bar{\rho}_W}{14\beta L\sqrt{\bar{\rho}_W}} \right\} .
\end{equation}
Then, we have, for all $k \geqslant 0$,
\begin{equation}
\begin{aligned}
&\mathbb{E} \left[ \left\| \tilde{\mathbf{y}}_{\beta \left( k+1 \right)} \right\| ^2 \right] 
\\
&\leqslant \frac{3+\bar{\rho}_W}{4}\mathbb{E} \left[ \left\| \tilde{\mathbf{y}}_{\beta k} \right\| ^2 \right] +\frac{18\bar{\rho}_WL^2}{1-\bar{\rho}_W}\mathbb{E} \left[ \left\| \tilde{\mathbf{x}}_{\beta k} \right\| ^2 \right] 
\\
&\quad +\frac{96n\gamma ^2\beta ^{2}L^2\bar{\rho}_W}{1-\bar{\rho}_W}\mathbb{E} \left[ \left\| \nabla f\left( \bar{x}_{\beta k} \right) \right\| ^2 \right] +6n\sigma ^2.
\end{aligned}
\end{equation}
\end{Lem}

\begin{proof}
Applying the recursion of FlexGT \eqref{Eq_y_update}, by Assumption~\ref{Ass_bounded_var}, we have
\begin{equation}
\begin{aligned}
&\mathbb{E} \left[ \left\| \tilde{\mathbf{y}}_{\beta \left( k+1 \right)} \right\| ^2 \right] 
\\
&\leqslant \frac{1+\bar{\rho} _{W}}{2}\mathbb{E} \left[ \left\| \tilde{\mathbf{y}}_{\beta k} \right\| ^2 \right] +2n\bar{\rho} _{W}\sigma ^2
\\
&\quad +\frac{\left( 1+\bar{\rho} _{W} \right) \bar{\rho} _{W}}{1-\bar{\rho} _{W}}\mathbb{E} \left[ \left\| \nabla F_{\beta \left( k+1 \right)}-\nabla F_{\beta k} \right\| ^2 \right] 
\\
&\quad +2\bar{\rho} _{W}\mathbb{E} \left[ \left< \tilde{\mathbf{y}}_{\beta k},\nabla F_{\beta k}-G_{\beta k} \right> \right] 
\\
&\quad +2\bar{\rho} _{W}\mathbb{E} \left[ \left< \nabla F_{\beta k}-G_{\beta k},\nabla F_{\beta \left( k+1 \right)}-\nabla F_{\beta k} \right> \right],
\end{aligned}
\end{equation}
where we have used Young's inequality and $\mathbb{E} \left[ G_{\beta k} \right]=\nabla F_{\beta k}$.
Then, noticing that
\begin{equation*}
\begin{aligned}
&2\bar{\rho}_W\mathbb{E} \left[ \left< \tilde{\mathbf{y}}_{\beta k},\nabla F_{\beta k}-G_{\beta k} \right> \right] 
\\
&=2\bar{\rho}_W\mathbb{E} \left[ \left< \left( \bar{W}-\mathbf{J} \right) G_{\beta k},\nabla F_{\beta k}-G_{\beta k} \right> \right] 
\leqslant 2\bar{\rho}_Wn\sigma ^2,
\end{aligned}
\end{equation*}
and
\begin{equation*}
\begin{aligned}
&2\bar{\rho}_W\mathbb{E} \left[ \left< \nabla F_{\beta \left( k+1 \right)}, \nabla F_{\beta k}-G_{\beta k} \right> \right] 
\\
&\leqslant \bar{\rho}_W\mathbb{E} \left[ \left\| \nabla F_{\beta \left( k+1 \right)}-\nabla F_{\beta k} \right\| ^2 \right] +\bar{\rho}_Wn\sigma ^2,
\end{aligned}
\end{equation*}
we further obtain
\begin{equation}
\begin{aligned}
&\mathbb{E} \left[ \left\| \tilde{\mathbf{y}}_{\beta \left( k+1 \right)} \right\| ^2 \right] 
\\
&\leqslant \frac{1+\bar{\rho}_W}{2}\mathbb{E} \left[ \left\| \tilde{\mathbf{y}}_{\beta k} \right\| ^2 \right] 
+\frac{6\bar{\rho}_WL^2}{1-\bar{\rho}_W}\mathbb{E} \left[ \left\| \tilde{\mathbf{x}}_{\beta k} \right\| ^2 \right] 
\\
&\quad +\frac{6\bar{\rho}_WL^2}{1-\bar{\rho}_W}\mathbb{E} \left[ \left\| \mathbf{x}_{\beta \left( k+1 \right)}-\mathbf{1}\bar{x}_{\beta k} \right\| ^2 \right] +5n\bar{\rho}_W\sigma ^2,
\end{aligned}
\end{equation}
where we have used the smoothness of the objectives.
Using Lemma~\ref{Lem_tech_client_diver} to bound the second last term and letting the stepsize satisfy \eqref{Eq_step_tracking_err}, we complete the proof.
\end{proof}

From Lemmas~\ref{Lem_cons_err} and \ref{Lem_tra_err}, it can be found that there is an interdependence between the consensus error and the GT error. To decouple these two errors, we introduce Lemmas~\ref{Lem_cons_err_sum_FlexGT} and \ref{Lem_trac_err_sum_FlexGT}, which bound the accumulated consensus and gradient tracking errors, respectively. For convenience, we initialize the consensus error to zero, i.e., $x_{i,0}=x_0 \in \mathbb{R} ^p, i\in \left[ n \right]$.

\begin{Lem}[\textbf{Accumulated consensus error}]\label{Lem_cons_err_sum_FlexGT}
    Suppose Assumptions~\ref{Ass_smooth}--\ref{Ass_graph} hold. Let the stepsize satisfy $\gamma \leqslant \frac{1}{4\beta L}$.
    Then, for all $K \geqslant 1$, we have
    \begin{equation}
    \begin{aligned}
&\frac{1}{K}\sum_{k=0}^{K-1}{\mathbb{E} \left[ \left\| \tilde{\mathbf{x}}_{\beta k} \right\| ^2 \right]}
\\
&\leqslant \frac{2\mathbb{E} \left[ \left\| \tilde{\mathbf{x}}_0 \right\| ^2 \right]}{\left( 1-\bar{\rho}_W \right) K}+\frac{16n\gamma ^2\beta ^{2}\bar{\rho}_W}{\left( 1-\bar{\rho}_W \right) ^2}\sigma ^2
\\
&\quad +\frac{8\gamma ^2\beta ^{2}\bar{\rho}_W}{\left( 1-\bar{\rho}_W \right) ^2K}\sum_{k=0}^{K-1}{\mathbb{E} \left[ \left\| \tilde{\mathbf{y}}_{\beta k} \right\| ^2 \right]}.
    \end{aligned}
    \end{equation}
\end{Lem}

\begin{proof}
Summing the obtained inequality in Lemma~\ref{Lem_cons_err} from $0$ to $K-1$, we get
\begin{equation}
 \begin{aligned}
&\sum_{k=0}^{K-1}{\mathbb{E} \left[ \left\| \tilde{\mathbf{x}}_{\beta \left( k+1 \right)} \right\| ^2 \right]}
\\
&\leqslant \frac{1+\bar{\rho}_W}{2}\sum_{k=0}^{K-1}{\mathbb{E} \left[ \left\| \tilde{\mathbf{x}}_{\beta k} \right\| ^2 \right]}+\frac{8n\gamma ^2\beta ^{2}\bar{\rho}_WK}{1-\bar{\rho}_W}\sigma ^2
\\
&\quad +\frac{4\gamma ^2\beta ^{2}\bar{\rho}_W}{1-\bar{\rho}_W}\sum_{k=0}^{K-1}{\mathbb{E} \left[ \left\| \tilde{\mathbf{y}}_{\beta k} \right\| ^2 \right]}  .
\end{aligned}
\end{equation}
Adding $\mathbb{E} \left[ \left\| \tilde{\mathbf{x}}_{0} \right\| ^2 \right]$ on both sides and noticing that
\[\sum_{k=0}^{K-1}{\mathbb{E} \left[ \left\| \tilde{\mathbf{x}}_{\beta k} \right\| ^2 \right]}\leqslant \sum_{k=0}^K{\mathbb{E} \left[ \left\| \tilde{\mathbf{x}}_{\beta k} \right\| ^2 \right]},
\]
we can obtain the target result, which completes the proof.
\end{proof}

\begin{Lem}[\textbf{Accumulated GT error}]\label{Lem_trac_err_sum_FlexGT}
    Suppose Assumptions~\ref{Ass_smooth}--\ref{Ass_graph} hold. Let the stepsize satisfy \eqref{Eq_step_tracking_err}.
    Then, for all $K \geqslant 1$, we have
\begin{equation}
\begin{aligned}
&\frac{1}{K}\sum_{k=0}^{K-1}{\mathbb{E} \left[ \left\| \tilde{\mathbf{y}}_{\beta k} \right\| ^2 \right]}
\\
&\leqslant \frac{4\mathbb{E} \left[ \left\| \tilde{\mathbf{y}}_0 \right\| ^2 \right]}{\left( 1-\bar{\rho}_W \right) K}+\frac{24n\sigma ^2}{\left( 1-\bar{\rho}_W \right) K}
\\
&\quad +\frac{96\bar{\rho}_WL^2}{\left( 1-\bar{\rho}_W \right) ^2}\frac{1}{K}\sum_{k=0}^{K-1}{\mathbb{E} \left[ \left\| \tilde{\mathbf{x}}_{\beta k} \right\| ^2 \right]}
\\
&\quad +\frac{384n\gamma ^2\beta ^{2}L^2\bar{\rho}_W}{\left( 1-\bar{\rho}_W \right) ^2}\frac{1}{K}\sum_{k=0}^{K-1}{\mathbb{E} \left[ \left\| \nabla f\left( \bar{x}_{\beta k} \right) \right\| ^2 \right]}.
\end{aligned}
\end{equation}
\end{Lem}

\begin{proof}
Based on Lemma~\ref{Lem_tra_err} and following a similar proof process as in Lemma~\ref{Lem_cons_err_sum_FlexGT}, summing the result from $0$ to $K-1$, we have
\begin{equation}
\begin{aligned}
&\sum_{k=0}^K{\mathbb{E} \left[ \left\| \tilde{\mathbf{y}}_{\beta k} \right\| ^2 \right]}
\\
&\leqslant \frac{3+\bar{\rho}_W}{4}\sum_{k=0}^{K-1}{\mathbb{E} \left[ \left\| \tilde{\mathbf{y}}_{\beta k} \right\| ^2 \right]}+\frac{18\bar{\rho}_WL^2}{1-\bar{\rho}_W}\sum_{k=0}^{K-1}{\mathbb{E} \left[ \left\| \tilde{\mathbf{x}}_{\beta k} \right\| ^2 \right]}
\\
&\quad +\frac{96n\gamma ^2\beta ^{2}L^2\bar{\rho}_W}{1-\bar{\rho}_W}\sum_{k=0}^{K-1}{\mathbb{E} \left[ \left\| \nabla f\left( \bar{x}_{\beta k} \right) \right\| ^2 \right]}+6Kn\sigma ^2.
\end{aligned}
\end{equation}
Adding $\mathbb{E} \left[ \left\| \tilde{\mathbf{y}}_{0} \right\| ^2 \right]$ on the both sides and noticing that
\[\sum_{k=0}^{K-1}{\mathbb{E} \left[ \left\| \tilde{\mathbf{y}}_{\beta k} \right\| ^2 \right]}\leqslant \sum_{k=0}^K{\mathbb{E} \left[ \left\| \tilde{\mathbf{y}}_{\beta k} \right\| ^2 \right]},
\]
we can obtain the desired result.
\end{proof}

\section*{References}

\bibliographystyle{ieeetr}
\bibliography{reference}

\end{document}